\documentclass[a4paper,11pt,twoside]{amsart}
\usepackage[english]{babel}
\usepackage[utf8]{inputenc}

\usepackage[a4paper,inner=3cm,outer=3cm,top=4cm,bottom=4cm,pdftex]{geometry}
\usepackage{fancyhdr}
\usepackage{color}
\usepackage{bold-extra}
\usepackage{ mathrsfs }
\usepackage{comment}
\usepackage{graphics}
\usepackage{aliascnt}
\usepackage[pdftex,citecolor=green,linkcolor=red]{hyperref}

\usepackage{amsmath}
\usepackage{amsfonts}
\usepackage{amssymb}
\usepackage{amsthm}
\usepackage{comment}
\usepackage{mathtools}
\usepackage{enumerate}
\usepackage{enumitem} 
\usepackage{bm}

\newtheorem{theorem}{Theorem}[section]

\newtheorem{lemma}[theorem]{Lemma}
\newtheorem{proposition}[theorem]{Proposition}
\theoremstyle{definition}
\newtheorem{definition}[theorem]{Definition}
\newtheorem{remark}[theorem]{Remark}

\numberwithin{equation}{section}

\newcommand\eps{\varepsilon}

\newcommand\E{\mathbb{E}}

\newcommand\R{\mathbb{R}}
\newcommand\Z{\mathbb{Z}}

\newcommand{\CC}{\mathcal{C}}

\newcommand{\CP}{\mathcal{P}}

\newcommand{\CF}{\mathcal{F}}

\newcommand{\CB}{\mathcal{B}}

\newcommand{\CD}{\mathcal{D}}
\newcommand{\CL}{\mathcal{L}}

\newcommand{\wt}{\widetilde}
\newcommand{\ve}{\mathbf}
\newcommand{\Mod}[1]{\ (\mathrm{mod}\ #1)}

\newcommand{\op}{\operatorname{op}}

\begin{document}
\title{Narrower arithmetic progressions in the primes}

\author[Shao]{Xuancheng Shao}
\address{Department of Mathematics, University of Kentucky\\
Lexington, KY 40506\\
USA}
\email{xuancheng.shao@uky.edu}
\thanks{XS was supported by NSF grant DMS-2452462.}

\begin{abstract}
For each fixed $k\geq 3$ and $\delta > 0$, we prove that every subset of the primes up to $X$ of relative density at least $\delta$ contains a non-trivial $k$-term arithmetic progression with common difference $O_{k,\delta,\eps}((\log X)^{k\cdot 2^{k-3}+\eps})$, for any $\eps>0$. This improves the previous exponent $(k-1)2^{k-2}$. The main ingredient is a relative Szemer\'{e}di theorem for narrow progressions under weaker linear forms conditions.
\end{abstract}

\maketitle
%\tableofcontents

\section{Introduction}

Let $\CP$ be the set of primes. Szemer\'{e}di's theorem \cite{Szemeredi1975} asserts that, for every fixed $k$ and $\delta>0$, every subset of $[X] := \{1,2,\cdots,X\}$ of size at least $\delta X$ contains a non-trivial $k$-term arithmetic progression ($k$-AP) once $X$ is sufficiently large in terms of $k,\delta$. Green and Tao \cite{GreenTao2008} extended this conclusion to subsets of $\CP\cap[X]$ of positive relative density. Their proof uses a relative Szemer\'{e}di theorem, which transfers the dense theorem to functions dominated by a suitable pseudorandom majorant. Conlon, Fox, and Zhao \cite{ConlonFoxZhao2015} subsequently weakened the pseudorandomness hypotheses, and Zhao \cite{Zhao2014} gave an arithmetic transference proof using a dense model theorem and a counting lemma. 

There are also sharp sparse random versions, due independently to Conlon and Gowers \cite{ConlonGowers2016} and Schacht \cite{Schacht2016}. More precisely, if $A\subset [X]$ is formed by retaining each integer independently with probability $p$, then for every fixed $k\geq3$ and $\delta\in(0,1)$ there exists $C=C(k,\delta)>0$ such that, asymptotically almost surely, every subset of $A$ of relative density at least $\delta$ contains a non-trivial $k$-AP whenever $p\geq C X^{-1/(k-1)}$. This threshold is sharp up to the constant $C$; see \cite[Theorem 1.12]{ConlonGowers2016} and \cite[Theorem 2.2]{Schacht2016}.

Quantitative versions of Szemer\'{e}di's theorem seek upper bounds for the largest size $r_k(X)$ of a subset of $[X]$ containing no non-trivial $k$-AP. Gowers \cite{Gowers2001} proved $r_k(X)\ll_k X(\log\log X)^{-c_k}$ for some $c_k>0$. Subsequent improvements for general $k$ include the bound $r_k(X)\ll_k X\exp(-(\log\log X)^{c_k})$ for some $c_k>0$ due to Leng, Sah, and Sawhney \cite{LengSahSawhney2024}.

In this paper we seek $k$-APs in dense subsets of $\CP \cap [X]$ with a small common difference $d \leq (\log X)^B$ for some constant $B$ in terms of $k$, loosely referred to as narrow $k$-APs. This problem is related to finding $k$-APs in dense subsets of primes in almost all short intervals $[x, x+(\log X)^B]$ with $x \sim X$. Heuristically, the proportion of primes in each such short interval is $\sim 1/(\log X)$, and hence the set of primes in such a short interval of length $M = (\log X)^B$ has density $M^{-1/B}$. Based on the random version of Szemer\'{e}di's theorem mentioned above, we expect the threshold $B > k-1$ which guarantees the existence of narrow $k$-APs. Tao and Ziegler \cite{TaoZiegler2015} proved the existence of a constant $B \ll_k 1$. The author \cite{Shao2017} subsequently obtained the exponent $B=(k-1)2^{k-2} + \eps$ for any fixed $\eps > 0$. In the current work, we obtain the exponent $B = k\cdot 2^{k-3} + \eps$.

\begin{theorem}[Narrow $k$-APs in primes]\label{thm:narrow-kap}
Fix a positive integer $k \geq 3$ and reals $\delta,\eps > 0$. Set $B_k = k \cdot 2^{k-3}$. If $X$ is sufficiently large in terms of $k,\delta,\eps$, then every subset $A \subset \CP\cap [X]$ with $|A| \geq \delta X/\log X$ contains a non-trivial $k$-term arithmetic progression $p,p+d,\cdots,p+(k-1)d$ with $1 \leq d\ll_{k,\delta,\eps}(\log X)^{B_k+\eps}$.
\end{theorem}

As in \cite{GreenTao2008,Shao2017}, we use a relative Szemer\'{e}di theorem, which transfers a positive lower bound for $k$-AP counts from bounded functions to functions dominated by pseudorandom majorants. Our improvement comes from weakening the linear forms conditions required for this transfer, so that the same majorant satisfies them at a shorter scale. We formulate the relative theorem on a prime cyclic group $G = \Z/N\Z$. 

\begin{theorem}[Relative Szemer\'{e}di theorem for narrow $k$-APs]\label{thm:relative-narrow-kap}
Fix positive integers $k,m\geq 3$ and a real $\delta > 0$. Let $N$ be a prime sufficiently large in terms of $k,m,\delta$, $G = \Z/N\Z$, and let $\eta > 0$ be sufficiently small in terms of $k,m,\delta$.
Let $\eta^{-2} \leq D \leq N$. If $\nu:G\rightarrow [0,+\infty)$ is a non-negative function satisfying the $(k,m)$-linear forms conditions of level $\eta$ at scale $\eta D$, then every function $f: G\rightarrow\R$ with $0 \leq f \leq \nu$ and $\E f\geq \delta$ satisfies
$$
\Lambda_D(f) := \E_{n \in G} \E_{d \in [D]} \prod_{j=1}^kf(n+(j-1)d) \geq c(k,\delta)
$$
for some constant $c(k,\delta)>0$.
\end{theorem}

The linear forms conditions appearing in the statement will be defined in Section \ref{sec:linear-forms}. Theorem \ref{thm:relative-narrow-kap} is in the exact same form as the relative Szemer\'{e}di theorems in \cite{ConlonFoxZhao2015,Zhao2014,Shao2017}, apart from the precise definitions of the linear forms conditions (and the extra parameter $m$ in our version). 

\subsection*{Related results}

Even for the full set of primes $A = \CP\cap [X]$, the range for the common difference $d$ in Theorem \ref{thm:narrow-kap} is already new. The Hardy--Littlewood prime tuples conjecture predicts infinitely many $k$-APs with a fixed common difference $d\ll_k1$. For $k=2$, bounded common differences follow from the work of Zhang \cite{Zhang2014} and Maynard \cite{Maynard2015}. Maynard also proved that bounded intervals can contain arbitrarily many primes. 

Gowers uniformity estimates also yield $k$-APs in primes in almost all short intervals. Matom\"aki, Radziwi\l\l, the author, Tao, and Ter\"av\"ainen \cite{MatomakiRadziwillShaoTaoTeravainen2026} prove higher-order uniformity of the primes in almost all intervals $(x,x+H]$, with $x\in[X,2X]$ and $X^{1/3+\eps}\leq H\leq X^{1-\eps}$. Combined with the generalized von Neumann theorem, these estimates give asymptotic counts of $k$-APs in primes in almost all such intervals. For all intervals, analogous asymptotic formulae were established when $X^{5/8+\eps}\leq H\leq X^{1-\eps}$ in \cite{MatomakiShaoTaoTeravainen2023}; the exponent $5/8$ has recently been improved to $3/5$ by Matom\"aki, Pandey, Pliego, Ter\"av\"ainen, and Wang \cite{MatomakiPandeyPliegoTeravainenWang2026}.

Tao and Ziegler's result \cite{TaoZiegler2015} also treats polynomial progressions $n+P_1(d),\cdots,n+P_k(d)$, where the integer polynomials $P_1,\cdots,P_k$ have a common constant term. They obtain $0<|d|\ll(\log X)^B$ for some constant $B$ depending on the polynomials and the relative density. This refines their earlier polynomial extension of the Green--Tao theorem \cite{TaoZiegler2008}, which gives such configurations with $d\leq X^{\eps}$ for any fixed $\eps>0$.

\subsection*{Organization of the paper}

In Section \ref{sec:outline} we give an outline of the proof of the relative Szemer\'{e}di theorem (Theorem \ref{thm:relative-narrow-kap}). This follows a standard transference argument, along with a new input (Proposition \ref{prop:new}) using our new linear forms conditions and a majorant for the $W$-tricked primes (Proposition \ref{prop:majorant}) satisfying these conditions. In Section \ref{sec:linear-forms} we specify the precise linear forms conditions, compare them with prior works, and give some intuition behind why they are weaker conditions. In Section \ref{sec:prime-majorant} we verify that the standard sieve majorant for $W$-tricked primes satisfies the required linear forms conditions, once the parameter $D$ exceeds its specified threshold. Via standard manipulations of the sieve majorant, this task is reduced to estimating a quantity ($L(\Psi)$ and its variant $L_H(\Psi)$) associated purely with the collection of linear forms involved. In Section \ref{sec:new}, we prove our new input (Proposition \ref{prop:new}), including an illustration of the idea in the special case $k=3$. In Section \ref{sec:transference}, we give the standard proofs of the remaining transference arguments. In Section \ref{sec:final-deduction}, we give the standard deduction of Theorems \ref{thm:narrow-kap} and \ref{thm:relative-narrow-kap}.

\subsection*{AI assistance}

ChatGPT (OpenAI) was used by the author to assist with developing mathematical arguments and preparing the manuscript. The author takes full responsibility for the results and their presentation.

\section{Proof outline of the relative Szemer\'{e}di theorem}\label{sec:outline}

\subsection{Notation and conventions}

We fix the positive integer $k \geq 3$ and work in $G = \Z/N\Z$, where $N$ is prime and tends to infinity. For a finite set, write $\E_{x \in X}$ for the average $|X|^{-1}\sum_{x\in X}$. When $X = G$, we sometimes omit $x \in G$ and simply write $\E$ for $\E_{x\in G}$. All functions on $G$ are real-valued. We use the probability measure on $G$ to define inner products and norms for functions on $G$:
$$
\langle f,g\rangle = \E_{n\in G}f(n)g(n), \ \ \|f\|_p^p = \E_{n\in G}|f(n)|^p.
$$
Following \cite{Shao2017}, we define the narrow $k$-AP counting operator at scale $D$ to be
$$
\Lambda_D(f_1,\cdots,f_k) = \E_{n \in G}\E_{d \in [D]}\prod_{j=1}^kf_j(n+(j-1)d).
$$
For $1 \leq i \leq k$ and a family $(u_j)_{1 \leq j \leq k}$ of functions on $G$, we define the dual function $\CD_i((u_j)_{j\neq i})$ on $G$ by
$$
\CD_i((u_j)_{j\neq i})(n) = \E_{d \in [D]}\prod_{j\neq i}u_j(n+(j-i)d).
$$
Note that the dependence of $\CD_i$ on the scale $D$ is suppressed for convenience. Thus
$$
\Lambda_D(f_1,\cdots,f_k) = \langle f_i, \CD_i((f_j)_{j\neq i})\rangle.
$$
We define the norm $\|\cdot\|_{D,i}$ of a function $h$ on $G$ by
$$
\|h\|_{D,i} = \sup_{|u_j| \leq 1} |\langle h,\CD_i((u_j)_{j\neq i})\rangle|.
$$
Thus if $f_1,\cdots,f_k$ are $1$-bounded and the function $f_i$ in the counting operator $\Lambda_D(f_1,\cdots,f_k)$ is replaced by a new function $g_i$ (while keeping the other functions unchanged), then the counting operator changes in value by at most $\|f_i-g_i\|_{D,i}$.

\subsection{The transference ingredients}

The transference argument has two parts. First the dense model lemma (Lemma \ref{lem:dense-model}) approximates a function $f$ dominated by a majorant $\nu$ by a $1$-bounded function $\wt{f}$ which is close to $f$ in the $\|\cdot\|_{D,i}$
norm for each $i$. This step requires $\nu$ to satisfy the ``mild extra conditions" in Definition \ref{def:linear-forms-conditions}. Second, the counting lemma (Lemma \ref{lem:counting}) shows that replacing $f$ by $\wt{f}$ preserves the counting operator, up to a negligible error. This step requires $\nu$ to satisfy two hypotheses which will be supplied by our new input (Proposition \ref{prop:new}) below.

\begin{lemma}[Dense model]\label{lem:dense-model}
Fix a positive integer $k \geq 3$ and a real $\eps > 0$. Let $N$ be a prime sufficiently large in terms of $k,\eps$, and let $\eta>0$ be sufficiently small in terms of $k,\eps$.
Let $\eta^{-2} \leq D \leq N$ and $1 \leq S \leq \eta D$. Let $\nu:G\rightarrow [0,+\infty)$ be a non-negative function satisfying
$$
\Big|\E_{n \in G}\E_{\ve{s}^{(0)},\ve{s}^{(1)}\in [S]^{k-1}} \prod_{\bm{\omega}\in\{0,1\}^{k-1}}\Big(\nu\Big(n+k!\sum_{j=1}^{k-1}s_j^{(\omega_j)}\Big)-1\Big)\Big| \leq \eta.
$$
Then for every function $f: G\rightarrow \R$ with $0 \leq f \leq \nu$ and $\E f \leq 1$, there exists $\wt{f}:G\rightarrow [0,1]$ with $\E\wt{f} = \E f$ such that 
$$
\|f - \wt{f}\|_{D,i} \leq \eps
$$
for each $1 \leq i \leq k$.
\end{lemma}

\begin{lemma}[Counting lemma by densification]\label{lem:counting}
Fix a positive integer $k\geq 3$. Let $N$ be a prime sufficiently large in terms of $k$, and let $\eta > 0$ be sufficiently small in terms of $k$.
Let $\eta^{-2} \leq D \leq N$. Let $\nu:G\rightarrow [0,+\infty)$ be a non-negative function satisfying the following two properties:
\begin{enumerate}
\item For $\nu_1,\cdots,\nu_k \in \{1,\nu\}$ and all $1 \leq i \leq k$, we have
$$
\|\CD_i((\nu_j)_{j\neq i})-1\|_2 \leq \eta.
$$
\item For all functions $f_j^{(\tau)}:G\rightarrow\R$ with $|f_j^{(\tau)}| \leq \nu+1$ and all $1 \leq i \leq k$, we have
$$
\Big|\E_{n \in G}(\nu(n)-1) \CD_i((f_j^{(0)})_{j\neq i})(n) \CD_i((f_j^{(1)})_{j\neq i})(n)\Big| \leq \eta.
$$
\end{enumerate}
For $1 \leq i \leq k$, let $f_i,\wt{f}_i,\nu_i:G\rightarrow\R$ be functions  such that $\nu_i\in\{1,\nu\}$, $0 \leq f_i \leq \nu_i$, and $0 \leq \wt{f}_i\leq 1$. If $\|f_i-\wt{f}_i\|_{D,i}\leq \eta$ for each $1\leq i \leq k$, then
$$
\Lambda_D(f_1,\cdots,f_k) = \Lambda_D(\wt{f}_1,\cdots,\wt{f}_k) + O_k(\eta^{1/2^k}).
$$
\end{lemma}

The precise statements of Lemmas \ref{lem:dense-model} and \ref{lem:counting} do not appear to be in the literature, but their proofs follow existing arguments in \cite{ConlonFoxZhao2015,Zhao2014,Shao2017}. We present their proofs in Section \ref{sec:transference} to highlight precisely what properties of $\nu$ are used in these two steps.

\subsection{The pseudorandom estimates and the prime majorant}

Our main new contribution is a verification of hypothesis (2) in Lemma \ref{lem:counting}, using a family of weaker linear forms conditions to be defined in Definition \ref{def:linear-forms-conditions}.

\begin{proposition}\label{prop:new}
Fix positive integers $k,m\geq 3$. Let $N$ be a prime sufficiently large in terms of $k,m$, and let $\eta > 0$ be sufficiently small in terms of $k,m$.
Let $\eta^{-2} \leq D \leq N$. Let $\nu:G\rightarrow [0,+\infty)$ be a non-negative function satisfying the $(k,m)$-linear forms conditions of level $\eta$ at scale $\eta D$. Then the following two conclusions hold:
\begin{enumerate}
\item For $\nu_1,\cdots,\nu_k \in \{1,\nu\}$ and all $1 \leq i \leq k$, we have
$$
\|\CD_i((\nu_j)_{j\neq i})-1\|_2 \ll_k \eta^{1/2}.
$$
\item For all functions $f_j^{(\tau)}:G\rightarrow\R$ with $|f_j^{(\tau)}| \leq \nu+1$ and all $1 \leq i \leq k$, we have
$$
\Big|\E_{n \in G}(\nu(n)-1) \CD_i((f_j^{(0)})_{j\neq i})(n) \CD_i((f_j^{(1)})_{j\neq i})(n)\Big| \ll_{k,m} \eta^{1/(m\cdot 2^{k-2})}.
$$
\end{enumerate}
\end{proposition}

We prove Proposition \ref{prop:new} in Section \ref{sec:new}. The heart of the matter is the second conclusion. The usual argument repeatedly applies Cauchy--Schwarz until all the functions $f_j^{(\tau)}$ have been eliminated; see \cite[Section 9]{Shao2017}. Here we stop after $k-3$ applications, leaving a bilinear average in two short variables. At this point, instead of bounding the bilinear average by the usual two applications of Cauchy--Schwarz, we bound it by an operator norm and use the inequality $\|T\|_{\op}^{2m}\leq\operatorname{tr}((TT^*)^m)$. This is the Schatten norm bound, previously used in \cite{Pascadi2026} in the related context of sparse regularity in graph theory (with arithmetic applications). The usual Cauchy--Schwarz approach corresponds to $m=2$, and we take higher trace moments.
Expanding the trace eventually produces the system $\CL_{m,i}$ in our linear forms conditions in Section \ref{sec:linear-forms}. 

The second new input is a majorant for the $W$-tricked primes satisfying our linear forms conditions.

\begin{proposition}\label{prop:majorant}
Fix positive integers $k,m\geq 3$, and let $\eta > 0$. Let $N$ be a prime sufficiently large in terms of $k,m,\eta$. Let $w \leq 0.1\log\log N$ be a slowly growing function of $N$, and let $W=\prod_{p\leq w}p$. Let $b\Mod{W}$ be a reduced residue class. There exists a non-negative function $\nu=\nu_{W,b}: G\rightarrow [0,+\infty)$ with the following properties.
\begin{enumerate}
\item (Majorant property) For $n \in G$ (viewed as an integer in $[N]$ in the natural way), if $Wn+b$ is prime and $Wn+b > N^{1/2}$, then
$$
\nu(n)\gg_{k,m} \frac{\varphi(W)}{W}\log N.
$$
\item (Linear forms conditions) The majorant satisfies the $(k,m)$-linear forms conditions of level $\eta$ at some scale 
$$
S \ll_{k,m,\eta} (\log N)^{k\cdot 2^{k-3}(1+ 4/m)}.
$$
\end{enumerate}
\end{proposition}

The construction is standard (see \cite{GreenTao2008,GreenTao2010,Shao2017}), and we follow the correlation estimates in \cite[Sections 3--5]{Shao2017} which eventually reduce the task to an analysis of the linear algebra properties of 
the linear systems involved. This is carried out in Section \ref{sec:prime-majorant}.

Once equipped with Lemma \ref{lem:dense-model}, Lemma \ref{lem:counting}, and Proposition \ref{prop:new}, the relative Szemer\'{e}di theorem (Theorem \ref{thm:relative-narrow-kap}) follows by combining them. Finally, Theorem \ref{thm:narrow-kap} follows by combining the relative Szemer\'{e}di theorem with the majorant from Proposition \ref{prop:majorant}, after choosing $m$ to be large enough. These deductions are given in Section \ref{sec:final-deduction}.

\section{The linear forms conditions}\label{sec:linear-forms}

In this section, we specify the linear systems used in our relative theorem.

\subsection*{The family $\Psi_i$ of linear forms}

For each $i \in [k]$, we fix once and for all two distinct indices $a=a_i,b=b_i\in [k]\setminus\{i\}$. Let $E = [k]\setminus\{i,a,b\}$. Let $\Psi_i$ be the collection of linear forms in the variables $u,v,\bm{\theta}=(d_0,d_1,(s_j^{(0)},s_j^{(1)})_{j\in E})$ given by
\begin{equation}\label{eq:psi_j}
\psi_j^{(\tau,\bm{\omega})}(\bm{\theta},u,v) = (j-i)d_\tau + (j-a)u + (j-b)v + \sum_{\ell\in E}(j-\ell)s_{\ell}^{(\omega_\ell)},
\end{equation}
where $j \in [k]$, $\tau\in\{0,1\}$, and $\bm{\omega}=(\omega_{\ell})_{\ell\in E}\in\{0,1\}^E$. We adopt the following conventions.
\begin{enumerate}
\item For $j=i$, $\psi_i^{(\tau,\bm{\omega})}$ does not depend on $\tau$, and so we simply write $\psi_i^{(\bm{\omega})}$. Thus there are $2^{k-3}$ linear forms in $\Psi_i$ of this form.
\item For $j\in E$, $\psi_j^{(\tau,\bm{\omega})}$ does not depend on $\omega_j$, and  so we may  impose $\bm{\omega}\in\{0,1\}^{E\setminus\{j\}}$. Thus there are $(k-3)\cdot 2\cdot 2^{k-4}$ linear forms in $\Psi_i$ of this form.
\item For $j\in\{a,b\}$, $\psi_a^{(\tau,\bm{\omega})}$ does not depend on $u$ and $\psi_b^{(\tau,\bm{\omega})}$ does not depend on $v$, and so we simply write $\psi_a^{(\tau,\bm{\omega})}(\bm{\theta},v)$ and $\psi_b^{(\tau,\bm{\omega})}(\bm{\theta},u)$. There are $2\cdot 2\cdot 2^{k-3}$ linear forms in $\Psi_i$ of this form.
\end{enumerate}
Thus in total we have
$$
|\Psi_i| = 2^{k-3} + (k-3)2^{k-3} + 4\cdot 2^{k-3} = (k+2)2^{k-3}.
$$

\subsection*{The set $\Theta_i$ of restricted parameters}

Let $\Theta_i$ be the set of parameters $\bm{\theta}$ such that all linear forms in $\Psi_i$ are distinct linear forms in $u,v$. Thus the complement of $\Theta_i$ is the union of proper hyperplanes of the form
$$
(j-i)d_{\tau} + \sum_{\ell\in E}(j-\ell)s_{\ell}^{(\omega_\ell)} = (j-i)d_{\tau'} + \sum_{\ell\in E}(j-\ell)s_{\ell}^{(\omega'_\ell)}
$$
for some $1 \leq j \leq k$, $\tau,\tau'\in\{0,1\}$, and $\bm{\omega},\bm{\omega}'\in\{0,1\}^E$, with $\psi_j^{(\tau,\bm{\omega})}\neq\psi_j^{(\tau',\bm{\omega}')}$.

\subsection*{The family $\CL_{m,i}$ of linear forms}

For each $i \in [k]$ and a positive integer $m\geq 3$, we use the same choice of $a,b$ as before and let $E = [k]\setminus\{i,a,b\}$ as before. Let $\CL_{m,i}$ be the collection of linear forms in the variables $\ve{u}=(u_1,\cdots,u_m)$, $\ve{v}=(v_1,\cdots,v_m)$, and $\bm{\theta}=(d_0,d_1,(s_j^{(0)},s_j^{(1)})_{j\in E})$ given by
\begin{equation}\label{eq:Lmi}
(\bm{\theta},\ve{u},\ve{v})\mapsto\psi_j^{(\tau,\bm{\omega})}(\bm{\theta},u_h,v_h) \text{ and } (\bm{\theta},\ve{u},\ve{v})\mapsto\psi_j^{(\tau,\bm{\omega})}(\bm{\theta},u_h,v_{h+1}),
\end{equation}
where $j \in [k]$, $\tau\in\{0,1\}$, and $\bm{\omega}\in\{0,1\}^E$,
with the convention that $v_{m+1}=v_1$. Each $\psi_j^{(\tau,\bm{\omega})} \in \Psi_i$ with $j\notin\{a,b\}$ leads to $2m$ distinct linear forms in $\CL_{m,i}$. However,  for $j=a$ there are only $m$ distinct linear forms $(\ve{\theta},\ve{u},\ve{v})\mapsto \psi_a^{(\tau,\bm{\omega})}(\ve{\theta},v_h)$ in $\CL_{m,i}$, and similarly for $j=b$. In total we have
$$
|\CL_{m,i}| = 2m(2^{k-3}+(k-3)2^{k-3}) + m\cdot 4\cdot 2^{k-3} = mk\cdot 2^{k-2}.
$$

\begin{definition}[Linear forms conditions]\label{def:linear-forms-conditions}
We use the notations described above. Fix positive integers $k,m\geq 3$. Let $N$ be a prime sufficiently large in terms of $k,m$. Let $2 \leq S \leq N$ and $\eta > 0$. We say that a non-negative function $\nu:G\rightarrow [0,+\infty)$ satisfies the $(k,m)$-linear forms conditions of level $\eta$ at scale $S$, if the following estimates hold simultaneously:
\begin{enumerate}
\item (Restricted conditions for $\CL_{m,i}$) For each $i \in [k]$, we have
$$
\Big|\E_{n,\bm{\theta},\ve{u},\ve{v}} 1_{\bm{\theta}\in\Theta_i} \prod_{\psi\in \CL_{m,i}} \nu(n+\psi(\bm{\theta},\ve{u},\ve{v})) - 1\Big| \leq \eta,
$$
where $n\in G$ and each variable in $\bm{\theta},\ve{u},\ve{v}$ runs over an interval of length at least $S$ contained in $[-2\eta^{-1}S, 2\eta^{-1}S]$. The same holds when $\psi$ runs over any subcollection of $\CL_{m,i}$.
\item (Conditions for $\Psi_i$) For each $i \in [k]$, we have
$$
\Big|\E_{n,\bm{\theta},u,v}\prod_{\psi\in \Psi_i} \nu(n+\psi(\bm{\theta},u,v)) - 1\Big| \leq \eta,
$$
where $n\in G$ and each variable in $\bm{\theta},u,v$ runs over an interval of length at least $S$ contained in $[-2\eta^{-1}S, 2\eta^{-1}S]$. The same holds when $\psi$ runs over any subcollection of $\Psi_i$.
\item (Mild extra conditions) Let $\CC$ be the collection of the following $2^{k-1}$ linear forms in $2(k-1)$ variables in $\ve{s}^{(0)} = (s_1^{(0)},\cdots,s_{k-1}^{(0)})$ and $\ve{s}^{(1)} = (s_1^{(1)},\cdots,s_{k-1}^{(1)})$:
$$
\CC = \Big\{k!\sum_{j=1}^{k-1}s_j^{(\omega_j)}: \bm{\omega} = (\omega_1,\cdots,\omega_{k-1})\in\{0,1\}^{k-1}\Big\}.
$$
Then we have
$$
\Big|\E_{n \in G}\E_{\ve{s}^{(0)},\ve{s}^{(1)}\in [S]^{k-1}} \prod_{\psi \in \CC}\nu(n+\psi(\ve{s}^{(0)},\ve{s}^{(1)})) - 1\Big| \leq \eta.
$$
The same holds when $\psi$ runs over any subcollection of $\CC$.
\end{enumerate}
\end{definition}

We give some intuition behind what our linear forms conditions look like and how they compare with \cite{Shao2017} in the special case $k=3$ and $i=1$ (with the choices $a=2$ and $b=3$). In this case $\Psi=\Psi_i$ consists of the five linear forms in the variables $u,v,d_0,d_1$ given by
\begin{equation}\label{eq:psik3}
-u-2v, \ \ d_0-v, \ \ d_1-v, \ \ 2d_0+u, \ \ 2d_1+u.
\end{equation}
For $\CL_m = \CL_{m,i}$, it has $6m$ linear forms in the variables $u_1,\cdots,u_m,v_1,\cdots,v_m,d_0,d_1$, which can be arranged on a cycle of length $2m$ as follows.

\begin{center}
\setlength{\unitlength}{1mm}\footnotesize
\begin{picture}(150,40)(0,-16)
\put(10,10){\line(1,0){38}}\put(48,10){\line(1,0){38}}
\put(102,10){\line(1,0){38}}
\put(10,10){\circle*{1.5}}\put(48,10){\circle*{1.5}}
\put(86,10){\circle*{1.5}}\put(102,10){\circle*{1.5}}\put(140,10){\circle*{1.5}}
\put(94,10){\makebox(0,0){$\cdots$}}
\qbezier(140,10)(140,-14)(75,-14)
\qbezier(75,-14)(10,-14)(10,10)
\put(6,8){\makebox(0,0){$v_1$}}\put(48,7){\makebox(0,0){$u_1$}}
\put(86,7){\makebox(0,0){$v_2$}}\put(102,7){\makebox(0,0){$v_m$}}
\put(144,8){\makebox(0,0){$u_m$}}
\put(10,18){\makebox(0,0){\shortstack{$d_0-v_1$\\$d_1-v_1$}}}
\put(48,18){\makebox(0,0){\shortstack{$2d_0+u_1$\\$2d_1+u_1$}}}
\put(86,18){\makebox(0,0){\shortstack{$d_0-v_2$\\$d_1-v_2$}}}
\put(102,18){\makebox(0,0){\shortstack{$d_0-v_m$\\$d_1-v_m$}}}
\put(140,18){\makebox(0,0){\shortstack{$2d_0+u_m$\\$2d_1+u_m$}}}
\put(29,12){\makebox(0,0){$-u_1-2v_1$}}
\put(67,12){\makebox(0,0){$-u_1-2v_2$}}
\put(121,12){\makebox(0,0){$-u_m-2v_m$}}
\put(75,-10){\makebox(0,0){$-u_m-2v_1$}}
\end{picture}
\end{center}

Thus the vertices on the cycle are $v_1,u_1,v_2,u_2,\cdots,v_m,u_m$. Each vertex is attached with two linear forms and each edge is attached with one linear form.
When $m=2$, there are $12$ linear forms in $\CL_2$:
$$
\renewcommand{\arraystretch}{1.4}
\begin{array}{cccc}
-u_1-2v_1 & -u_1-2v_2 & -u_2-2v_2 & -u_2-2v_1 \\
 d_0-v_1 & d_1-v_2 & d_0-v_2 & d_1-v_1 \\
 2d_0+u_1 & 2d_1+u_1 & 2d_0+u_2 & 2d_1+u_2
\end{array}
$$
Note that each column forms a $3$-AP. These are  precisely the $12$ linear forms in \cite[Example 2.2]{Shao2017} (after renaming the variables). Thus the main innovation in our work is the extension to larger $m$, along with the restriction to $\Theta = \Theta_i$. As $m$ grows,  even though there are more linear forms in $\CL_m$, the restricted conditions for $\CL_m$ in the $(k,m)$-linear forms conditions are actually easier to satisfy, as we will prove in Section \ref{sec:prime-majorant}.

\section{The pseudorandom majorant for $W$-tricked primes}\label{sec:prime-majorant}

In this section we prove Proposition \ref{prop:majorant}. The construction of $\nu=\nu_{W,b}$ is fairly standard as follows (see \cite[Appendix D]{GreenTao2010} and \cite[Section 3]{Shao2017}). Set $R = N^{1/(100^km)}$ (say). Let $\chi:\R\rightarrow [0,+\infty)$ be a fixed smooth function supported on $[-1,1]$, with $\chi(0) > 0$ and $\int_0^{\infty}|\chi'(x)|^2dx=1$. Define
$$
\Lambda_{\chi,R}(n) = (\log R)\sum_{d\mid n}\mu(d)\chi\Big(\frac{\log d}{\log R}\Big)
$$
and define the majorant by
$$
\nu(n) = \nu_{W,b}(n) = \frac{\varphi(W)}{W\log R}\Lambda_{\chi,R}(Wn+b)^2
$$
for $n \in [N]$. By the arguments in \cite{Shao2017}, we have $\E\nu = 1+o(1)$ and the majorant property holds. It suffices to prove for any $\eta > 0$ that $\nu$ satisfies the $(k,m)$-linear forms conditions of level $\eta$ at some scale
$$
S \ll_{k,m,\eta} (\log N)^{k\cdot 2^{k-3}(1+4/m)}.
$$

\subsection{The linear forms conditions for $\Psi_i$ and $\CC$}\label{sec:linear-forms-Psi}

We start with verifying (2) and (3) in Definition \ref{def:linear-forms-conditions}.
Fix $i \in [k]$ and write $\Psi = \Psi_i$ throughout this subsection. Recall that $\Psi$ is a collection of $(k+2)2^{k-3}$ linear forms in $2k-2$ variables given in \eqref{eq:psi_j}.
Recall the constant $L(\Psi)$ from \cite[Definition 2.4]{Shao2017} defined as follows. For a partition $\pi$ of $\Psi$ into $s = s(\pi)$ disjoint nonempty subsets (called blocks) $\Psi_1,\cdots,\Psi_s$, let $\Pi(\Psi,\pi)$ be the subspace on which linear forms in the same block $\Psi_i$ all agree:
$$
\Pi(\Psi,\pi) = \{\ve{x}\in\R^{2k-2}: \psi(\ve{x})=\psi'(\ve{x})\text{ whenever }\psi,\psi'\in \Psi_i\text{ for some }1 \leq i \leq s\}.
$$
Thus the codimension of $\Pi(\Psi,\pi)$ is equal to $\dim U_\pi$, where
$$
U_\pi = \operatorname{Span}\{\psi-\psi': \psi,\psi'\in\Psi_i\text{ for some }1 \leq i \leq s\}.
$$
Note that $\dim U_\pi=0$ if and only if $\pi$ is the trivial partition into singletons.
The constant $L(\Psi)$ is defined by
$$
L(\Psi) = \max_{\pi} \frac{|\Psi|-s(\pi)}{\dim U_\pi},
$$
where the maximum is taken over all nontrivial partitions $\pi$ (with the convention that an empty maximum is defined to be $0$). Note that $L(\Psi') \leq L(\Psi)$ for any subcollection $\Psi'\subset\Psi$. Indeed, if $\pi'$ is a partition of $\Psi'$, we may extend it by singletons to a partition $\pi$ of $\Psi$. Then $|\Psi|-s(\pi)=|\Psi'|-s(\pi')$ and $U_\pi = U_{\pi'}$.
Hence, in view of \cite[Proposition 2.6]{Shao2017}, it suffices to show that
$$
L(\Psi) \leq k\cdot 2^{k-3}.
$$
Let $\pi$ be a nontrivial partition of $\Psi$. It suffices to show that
$$
|\Psi|-s \leq k \cdot 2^{k-3}\cdot \dim U_\pi.
$$
If $\dim U_\pi \geq 2$, then this holds trivially since $|\Psi| = (k+2)2^{k-3}$. Henceforth assume that $\dim U_\pi=1$. Recall that the variables of the linear forms are $u,v,\bm{\theta}=(d_0,d_1,(s_j^{(0)},s_j^{(1)})_{j\in E})$. Consider
$$
V = \operatorname{Span}\Big(\{d_1-d_0\} \cup\{s_j^{(1)}-s_j^{(0)}: j \in E\}\Big).
$$

We divide into two cases. If $U_\pi\not\subseteq V$, then $U_\pi\cap V=\{0\}$. Hence the $2^{k-2}$ linear forms $\psi_a^{(\tau,\bm{\omega})}$ for $\tau\in\{0,1\}$ and $\bm{\omega}\in\{0,1\}^E$, which are of the form
$$
\psi_a^{(\tau,\bm{\omega})}(\bm{\theta},u,v) = (a-i)d_\tau+(a-b)v+\sum_{\ell\in E}(a-\ell)s_\ell^{(\omega_\ell)},
$$
remain distinct modulo $U_\pi$, and thus they must lie in different blocks of the partition $\pi$. It follows that $s \geq 2^{k-2}$, and hence
$$
|\Psi|-s \leq (k+2)2^{k-3} - 2^{k-2} = k\cdot 2^{k-3}.
$$

Now assume that $U_\pi\subset V$. Choose a nonzero element in $U_\pi$ of the form
$$
\lambda(d_0-d_1) + \sum_{\ell \in E}\lambda_\ell(s_\ell^{(1)}-s_\ell^{(0)})
$$
for some scalars $\lambda$ and $\lambda_\ell$ for $\ell \in E$. Thus if $\psi_j^{(\tau',\bm{\omega}')}-\psi_{j'}^{(\tau,\bm{\omega})} \in U_\pi$ for two distinct forms $\psi_j^{(\tau',\bm{\omega}')}$ and $\psi_{j'}^{(\tau,\bm{\omega})}$, then 
$j$ must equal $j'$, $\tau=\tau'$ (or $j=i$) if and only if $\lambda=0$, and $\omega_\ell=\omega_\ell'$ (or $\ell=j$) if and only if $\lambda_\ell=0$ for each $\ell \in E$. It follows that for each linear form $\psi\in\Psi$, there is at most one other form $\psi'\in\Psi$ with $\psi'-\psi\in U_\pi$. Hence each block of the partition $\pi$ has size at most $2$, and
$$
|\Psi|-s \leq \frac{1}{2}|\Psi| \leq \frac{1}{2}(k+2)2^{k-3} < k\cdot 2^{k-3}.
$$
This completes the proof of the linear forms conditions for $\Psi$. The mild extra conditions for $\CC$ also follow since one can easily deduce that $L(\CC) = 2^{k-2} < k\cdot 2^{k-3}$.

\begin{remark}\label{rem:L-constant1}
In fact, the exact value is
$$
L(\Psi)=(k+1)2^{k-4}. 
$$
We do not pursue this since the weaker bound $L(\Psi) \leq k\cdot 2^{k-3}$ obtained above suffices for our purposes. Here is the partition $\pi$ showing that $L(\Psi) \geq (k+1)2^{k-4}$. For each $j \neq i$, pair the forms $\psi_j^{(0,\bm{\omega})}$ and $\psi_j^{(1,\bm{\omega})}$ in the same block. Keep the forms $\psi_i^{(\bm{\omega})}$ as singleton blocks. Then one can verify that $|\Psi|-s(\pi) = (k+1)2^{k-4}$. Moreover, $U_\pi$ is the span of $d_1-d_0$, so that $\dim U_\pi=1$. 
\end{remark}

\subsection{The restricted linear forms conditions for $\CL_{m,i}$}

Now we verify (1) in Definition \ref{def:linear-forms-conditions}. First we record an analogue of \cite[Proposition 2.6]{Shao2017} when there are restrictions on parameters in the linear forms condition.

\begin{lemma}[Restricted linear forms conditions]\label{lem:restricted-linear-forms}
Let $\nu=\nu_{W,b}$ be constructed as above.
Fix a system $\Psi$ of at most $50^km$ distinct homogeneous linear forms in $d$ variables with integral coefficients, and
let $H\subset\R^d$ be a fixed finite union of proper hyperplanes.
Define
$$
L_H(\Psi)=\max_{\pi:\Pi(\Psi,\pi)\not\subseteq H}
\frac{|\Psi|-s(\pi)}{\dim U_\pi},
$$
where $s(\pi),U_\pi,\Pi(\Psi,\pi)$ are defined as in Section \ref{sec:linear-forms-Psi}, and the maximum is taken over all nontrivial partitions $\pi$ of $\Psi$ with $\Pi(\Psi,\pi)\not\subseteq H$ (with the convention that an empty maximum is defined to be $0$).
Let $N$ be a prime sufficiently large. For any $\eta > 0$, there exists
$$
S \ll_{\Psi,H,\eta} (\log N)^{L_H(\Psi)}
$$
such that the majorant $\nu=\nu_{W,b}$ constructed above satisfies
$$
\Big|\E_{n \in G} \E_{\ve{x} \in I_1\times\cdots\times I_d}
1_{\ve{x}\notin H}\prod_{\psi\in\Psi}\nu(n+\psi(\ve{x})) - 1\Big| \leq \eta,
$$
whenever $I_1,\cdots,I_d$ are intervals of length at least $S$ contained in $[-2\eta^{-1}S, 2\eta^{-1}S]$.
\end{lemma}

\begin{proof}
This is an extension of \cite[Proposition 2.6]{Shao2017}, which corresponds to the special case when $H$ is empty. The proof follows the argument in \cite[Section 5]{Shao2017}, and we only give the modifications needed. As in the beginning of \cite[Section 5]{Shao2017}, we split the average over $\ve{x}$ according to the equalities among the values of $\psi(\ve{x})$ for $\psi\in\Psi$. In particular, each $\ve{x}\in\R^d$ induces a partition $\pi(\ve{x})$ of $\Psi$, where two linear forms $\psi,\psi'\in\Psi$ lie in the same block if and only if $\psi(\ve{x}) = \psi'(\ve{x})$.
It suffices to show that for each partition $\pi$ we have
\begin{equation}\label{eq:restricted-sum}
\E_{n \in G} \E_{\ve{x}\in I_1\times\cdots\times I_d} 1_{\pi(\ve{x})=\pi} 1_{\ve{x}\notin H}\prod_{\psi\in\Psi}\nu(n+\psi(\ve{x})) = 1_{s(\pi)=|\Psi|} + O(\eta).
\end{equation}
If $\Pi(\Psi,\pi) \subseteq H$, then $\pi(\ve{x})=\pi$ implies $\ve{x}\in H$, and the left-hand side above vanishes. Thus we may assume that $\Pi(\Psi,\pi)\not\subseteq H$.

First consider the nontrivial partitions $\pi$, so that $s(\pi)<|\Psi|$. By non-negativity, it suffices to prove \eqref{eq:restricted-sum} without the restriction $\ve{x}\notin H$. The arguments in \cite{Shao2017} (see the final displayed equation in \cite[Section 5]{Shao2017}) therefore give an upper bound of
$$
\ll S^{-\dim U_\pi}(\log N)^{|\Psi|-s(\pi)}
$$
for the left-hand side of \eqref{eq:restricted-sum}. This can be made $O(\eta)$ once 
$$
S \geq \eta^{-1}(\log N)^{(|\Psi|-s(\pi))/\dim U_\pi}.
$$

It remains to consider the trivial partition $\pi$ of $\Psi$ into singletons. Let $X_0$ be the set of $\ve{x} \in \R^d$ such that $\pi(\ve{x})$ is trivial. The arguments in \cite{Shao2017} (see the penultimate displayed equation in \cite[Section 5]{Shao2017}) give
$$
\E_{n\in G}\E_{\ve{x} \in I_1\times\cdots\times I_d} 1_{\ve{x}\in X_0}\prod_{\psi\in\Psi}\nu(n+\psi(\ve{x})) = 1+o(1) + O(S^{-1}).
$$
Thus it suffices to show that
\begin{equation}\label{eq:restricted-H}
\Big|\E_{n \in G}\E_{\ve{x} \in I_1\times\cdots\times I_d}1_{\ve{x}\in X_0\cap H}\prod_{\psi\in\Psi}\nu(n+\psi(\ve{x}))\Big| = O(\eta).
\end{equation}
For fixed $\ve{x}\in X_0$, the values of $\psi(\ve{x})$ for $\psi \in \Psi$ are all distinct. Hence \cite[Proposition 3.3]{Shao2017} gives
$$
\Big|\E_{n \in G}\prod_{\psi\in\Psi}\nu(n+\psi(\ve{x}))\Big| \ll E(\Psi(\ve{x}))+o(1),
$$
where 
$$
E(\Psi(\ve{x})) = \exp\Big(O\Big(\sum_{p\mid\Delta(\Psi(\ve{x}))}\frac{1}{p}\Big)\Big),\ \ \Delta(\Psi(\ve{x})) = \prod_{\substack{\psi,\psi'\in\Psi \\ \psi(\ve{x})\neq\psi'(\ve{x})}}(\psi(\ve{x})-\psi'(\ve{x})).
$$
Squaring this bound we get
$$
\E_{\ve{x}} 1_{\ve{x}\in X_0}\Big|\E_{n \in G}\prod_{\psi\in\Psi}\nu(n+\psi(\ve{x}))\Big|^2 \ll \E_{\ve{x}} E(\Psi(\ve{x}))^2 + o(1),
$$
which is $O(1)$ by the estimates near the end of \cite[Section 5]{Shao2017}. By Cauchy--Schwarz, the square of the left-hand side of \eqref{eq:restricted-H} is bounded by
$$
\Big(\E_{\ve{x}} 1_{\ve{x}\in H}\Big) \Big(\E_{\ve{x}} 1_{\ve{x}\in X_0}\Big|\E_{n \in G}\prod_{\psi\in\Psi}\nu(n+\psi(\ve{x}))\Big|^2\Big) \ll \frac{|H \cap (I_1\times\cdots\times I_d)|}{|I_1|\cdots |I_d|}.
$$
Since $H$ is a union of fixed hyperplanes, the quantity above is $O(S^{-1})$. This concludes the proof of \eqref{eq:restricted-H} since we may choose $S$ sufficiently large in terms of $\eta$.
\end{proof}

Fix $i \in [k]$ and write $\Psi = \CL_{m,i}, H=\Theta_i^c$ throughout this subsection. Recall that $\Psi$ is a collection of $mk\cdot 2^{k-2}$ linear forms in $2k-4+2m$ variables given in \eqref{eq:Lmi}. Note that $L_H(\Psi') \leq L_H(\Psi)$ for any subcollection $\Psi'\subset\Psi$. In view of Lemma \ref{lem:restricted-linear-forms}, it suffices to prove that
$$
L_H(\Psi) \leq k\cdot 2^{k-3}\Big(1+\frac{4}{m}\Big). 
$$
Recall that variables in the linear forms are $\ve{u}=(u_1,\cdots,u_m)$, $\ve{v}=(v_1,\cdots,v_m)$, $\bm{\theta}=(d_0,d_1,(s_j^{(0)},s_j^{(1)})_{j\in E})$, and $H$ is the union of proper hyperplanes of the form
$$
(j-i)d_{\tau} + \sum_{\ell\in E}(j-\ell)s_{\ell}^{(\omega_\ell)} = (j-i)d_{\tau'} + \sum_{\ell\in E}(j-\ell)s_{\ell}^{(\omega'_\ell)}
$$
Let $\pi$ be a nontrivial partition of $\Psi$ with $\Pi(\Psi,\pi)\not\subseteq H$. Let $V$ be defined as in Section \ref{sec:linear-forms-Psi}:
$$
V = \operatorname{Span}\Big(\{d_1-d_0\} \cup\{s_j^{(1)}-s_j^{(0)}: j \in E\}\Big).
$$
We reduce the problem modulo $V$. Each $\psi_j^{(\tau,\bm{\omega})}\Mod{V}$ can be viewed as a linear form $\overline{\psi}_j$ in the variables $d,u,v,(s_j)_{j\in E}$ given by the formula
\begin{equation}\label{eq:psibar}
\overline{\psi}_j(u,v) = (j-i)d + (j-a)u + (j-b)v + \sum_{\ell\in E}(j-\ell)s_\ell,
\end{equation}
where we suppress the dependence on the other variables. For $j=a$ we simply write $\overline{\psi}_a(v)$, and for $j=b$ we simply write $\overline{\psi}_b(u)$. Let $\overline{\Psi}$ be the collection of distinct linear forms $\psi\Mod{V}$ with $\psi\in\Psi$. Thus
$$
\overline{\Psi} = \{\overline{\psi}_j(u_h,v_h), \overline{\psi}_j(u_h,v_{h+1}): j\notin\{a,b\},1 \leq h \leq m\} \cup \{\overline{\psi}_a(v_h),\overline{\psi}_b(u_h): 1 \leq h \leq m\}. 
$$
All linear forms in the fibre over $\overline{\psi}\in\overline{\Psi}$ (under the natural projection map $\Psi\rightarrow\overline{\Psi})$ are distinct modulo $U_\pi$. This follows from $\Pi(\Psi,\pi)\not\subseteq H$: the difference of two distinct forms in one fibre defines one of the hyperplanes of $H$, so this difference cannot belong to $U_\pi$. Note that for each $\overline{\psi}_a(v_h)$ and $\overline{\psi}_b(u_h)$, there are $2^{k-2}$ linear forms in the fibre, and for each $\overline{\psi}_j(u_h,v_h)$ and $\overline{\psi}_j(u_h,v_{h+1})$ with $j\notin\{a,b\}$, there are $2^{k-3}$ elements in the fibre.

We arrange the linear forms in $\overline{\Psi}$ in the following manner. On a cycle of length $2m$ with vertices labeled by the variables $v_1,u_1,v_2,u_2,\cdots,v_m,u_m$, we attach to each edge $(v_h,u_h)$ the $k-2$ linear forms $\overline{\psi}_j(u_h,v_h)$ for $j\notin\{a,b\}$, attach to each edge $(u_h,v_{h+1})$ the $k-2$ linear forms $\overline{\psi}_j(u_h,v_{h+1})$ for $j\notin\{a,b\}$, attach to each vertex $u_h$ the linear forms $\overline{\psi}_b(u_h)$, and attach to each vertex $v_h$ the linear forms $\overline{\psi}_a(v_h)$. For the rest of this proof, we rename the linear forms on the vertices $v_1,u_1,\cdots,v_m,u_m$ by $z_1,z_2,\cdots,z_{2m}$ in this order (with the convention that $z_{2m+1}=z_1$), and we also refer to the vertices by $z_1,\cdots,z_{2m}$. This is illustrated in the picture below.

\begin{center}
\setlength{\unitlength}{1mm}\footnotesize
\begin{picture}(150,44)(0,-22)
\put(10,10){\line(1,0){38}}\put(48,10){\line(1,0){38}}
\put(102,10){\line(1,0){38}}
\put(10,10){\circle*{1.5}}\put(48,10){\circle*{1.5}}
\put(86,10){\circle*{1.5}}\put(102,10){\circle*{1.5}}\put(140,10){\circle*{1.5}}
\put(94,10){\makebox(0,0){$\cdots$}}
\qbezier(140,10)(140,-12)(75,-12)
\qbezier(75,-12)(10,-12)(10,10)
\put(6,8){\makebox(0,0){$v_1$}}\put(48,7){\makebox(0,0){$u_1$}}
\put(86,7){\makebox(0,0){$v_2$}}\put(102,7){\makebox(0,0){$v_m$}}
\put(144,8){\makebox(0,0){$u_m$}}
\scriptsize
\put(10,17){\makebox(0,0){\shortstack{$z_1=\overline{\psi}_a(v_1)$}}}
\put(48,17){\makebox(0,0){\shortstack{$z_2=\overline{\psi}_b(u_1)$}}}
\put(86,17){\makebox(0,0){\shortstack{$z_3=\overline{\psi}_a(v_2)$}}}
%\put(102,19){\makebox(0,0){}
\put(140,17){\makebox(0,0){\shortstack{$z_{2m}=\overline{\psi}_b(u_m)$}}}
\put(29,13){\makebox(0,0){$w_{2,j}=\overline{\psi}_j(u_1,v_1)$}}
\put(67,13){\makebox(0,0){$w_{3,j}=\overline{\psi}_j(u_1,v_2)$}}
\put(121,13){\makebox(0,0){$w_{2m,j}=\overline{\psi}_j(u_m,v_m)$}}
\put(75,-8){\makebox(0,0){$w_{1,j}=\overline{\psi}_j(u_m,v_1)$}}
\put(75,-17){\makebox(0,0){$j\notin\{a,b\}$ on every edge}}
\end{picture}
\end{center}

Let $\overline{W}$ be the span of $z_1,\cdots,z_{2m}$.
Note that these $2m$ forms are linearly independent. From \eqref{eq:psibar} we see that the $k-2$ forms on an edge $(z_{h-1},z_h)$ are
$$
w_{h,j} := \frac{j-b}{a-b}z_h + \frac{j-a}{b-a}z_{h-1} \ \ \text{or} \ \ w_{h,j} := \frac{j-b}{a-b}z_{h-1} + \frac{j-a}{b-a}z_{h}
$$
for $j\notin\{a,b\}$.
Let $\overline{U}_\pi$ be the quotient $U_\pi/(V \cap U_\pi)$. We produce a basis of $\overline{W}/(\overline{W}\cap\overline{U}_\pi)$ by visiting the vertices $z_1,\cdots,z_{2m}$ in this order and adding $z_h + (\overline{W}\cap\overline{U}_\pi)$ to the basis if it is not in the span of the previous basis elements. Let $\CB\subset [2m]$ be the set of indices $h$ such that $z_h + (\overline{W}\cap\overline{U}_\pi)$ is included in this basis. Since $z_1\notin \overline{U}_\pi$ (as every generator $\overline{\psi}-\overline{\psi}'$ of $\overline{U}_\pi$ is a linear combination of $z_1,\cdots,z_h$ with coefficients summing to $0$), we have $1 \in \CB$ by construction. We have
$$
|\CB| = 2m - \dim(\overline{W}\cap\overline{U}_\pi) \geq 2m-\dim U_\pi.
$$
Define a subset $\overline{S}\subset\overline{\Psi}$ by including the forms $z_h$ in $\overline{S}$ for each $h \in \CB$, as well as the forms $w_{h,j}$ ($j\notin\{a,b\}$) on the edge $(z_{h-1},z_h)$ for each $h \in \CB\setminus\{1\}$. Thus
$$
\overline{S} = \bigcup_{h\in\CB}\overline{S}_h, \ \ \overline{S}_1=\{z_1\}, \ \ \overline{S}_h = \{z_h\}\cup\{w_{h,j}: j \notin\{a,b\}\}\text{ for }h\in\CB\setminus\{1\}.
$$
Let $S \subset \Psi$ be the fibre over $\overline{S}$. Then
\begin{equation}\label{eq:S-bound}
|S| = 2^{k-2}|\CB| + (k-2)2^{k-3}(|\CB|-1) \geq [k(2m-\dim U_\pi)-(k-2)]2^{k-3}.
\end{equation}
We claim that all forms in $S$ are distinct modulo $U_\pi$. Indeed, it suffices to show that all forms in $\overline{S}$ are distinct modulo $\overline{U}_\pi$, since each fibre contains distinct forms modulo $U_\pi$. Now pick two distinct forms, one from $\overline{S}_h$ and the other from $\overline{S}_{h'}$ for some $h,h' \in \CB$ with $h' \leq h$. If they are equal modulo $\overline{U}_\pi$, then this leads to a linear relation of the form
$$
(\lambda_h z_h + \lambda_{h-1}z_{h-1}) - (\lambda'_{h'}z_{h'} +\lambda'_{h'-1}z_{h'-1}) \in \overline{U}_\pi
$$
for some scalars $\lambda_h,\lambda_{h-1},\lambda'_{h'},\lambda'_{h'-1}$ with $\lambda_h,\lambda'_{h'} \neq 0$ (with the convention that $z_0=0$). If $h' < h$, this implies that $z_h\in\operatorname{Span}\{z_1,\cdots,z_{h-1}\}$ modulo $\overline{U}_\pi$, contradicting the construction of $\CB$. If $h'=h$, then we also have $\lambda_h\neq\lambda'_{h'}$, which still implies that $z_h\in\operatorname{Span}\{z_1,\cdots,z_{h-1}\}$ modulo $\overline{U}_\pi$, contradicting the construction of $\CB$.

This shows that all forms in $S$ are distinct modulo $U_\pi$, and hence they must lie in different blocks of the partition $\pi$, so that $s(\pi) \geq |S|$. Combining this with \eqref{eq:S-bound}, we have
$$
\frac{|\Psi|-s(\pi)}{\dim U_\pi} \leq \frac{[k(\dim U_\pi)+k-2]2^{k-3}}{\dim U_\pi} < k\cdot 2^{k-3}\Big(1 + \frac{1}{\dim U_\pi}\Big).
$$
This already gives the desired bound if $\dim U_\pi \geq m/4$.
Henceforth assume that $\dim U_\pi<m/4$. By the definition of $U_\pi$, we can pick a basis of $\overline{U}_\pi$ consisting of $<m/4$ forms of the form $\overline{\psi}-\overline{\psi}'$ with $\overline{\psi},\overline{\psi}'\in\overline{\Psi}$. Since each form in $\overline{\Psi}$ involves at most two vertices, all these basis elements involve $<m$ vertices. Hence there must be two adjacent vertices which do not occur in $\overline{U}_\pi$. By relabeling the vertices on the cycle if necessary, we may assume that $u_m,v_1$ do not occur in $\overline{U}_\pi$ (i.e. these two variables have zero coefficients in every form in $\overline{U}_\pi$.)

If $\overline{S}\subset\overline{\Psi}$ is the set of distinct linear forms modulo $\overline{U}_\pi$ constructed above, then $\overline{S}'=\overline{S}\cup\{\overline{\psi}_j(u_m,v_1): j \notin\{a,b\}\}$ still has this property, since any collision modulo $\overline{U}_\pi$ leads to a form in $\overline{U}_\pi$ with a nonzero coefficient at either $u_m$ or $v_1$, a contradiction. Let $S'$ be the fibre over $\overline{S}'$, so that the linear forms in $S'$ are distinct modulo $U_\pi$, and
$$
|S'| \geq |S| + (k-2)2^{k-3} \geq k(2m-\dim U_\pi)2^{k-3}
$$
by \eqref{eq:S-bound}. Since $s(\pi) \geq |S'|$, we have
$$
\frac{|\Psi|-s(\pi)}{\dim U_\pi} \leq \frac{k(\dim U_\pi)2^{k-3}}{\dim U_\pi} = k\cdot 2^{k-3}.
$$
This completes the proof.

\begin{remark}
The restriction $\ve\theta\in\Theta_i$ in the linear forms conditions for $\CL_{m,i}$ is crucial. In fact,
$$
L(\CL_{m,i})=m(k-1)2^{k-3},
$$
which increases linearly with $m$. The lower bound can be seen by pairing the $\tau=0,1$ forms for every $j\neq i$ and keeping the $j=i$ forms singleton (as in Remark \ref{rem:L-constant1}). The resulting partition $\pi$ has $U_\pi=\operatorname{Span}(d_1-d_0)$ and $|\CL_{m,i}|-s(\pi)=m(k-1)2^{k-3}$. When $m=2$, this matches the exponent $(k-1)2^{k-2}$ from \cite{Shao2017}.

On the other hand, for the restricted version $L_H(\CL_{m,i})$ (with $H = \Theta_i^c$), consider the partition $\pi$ obtained by putting all copies from each fixed $\psi\in\Psi_i$ into one block. Then $s(\pi)=|\Psi_i|=(k+2)2^{k-3}$ and
$$
U_{\pi}=\operatorname{Span}\{u_h-u_1,v_h-v_1:2\leq h\leq m\}.
$$
Thus $\dim U_{\pi}=2m-2$. Moreover, $\Pi(\CL_{m,i},\pi)$ is the subspace defined by $u_1=\cdots=u_m$ and $v_1=\cdots=v_m$. These equations do not involve $\bm{\theta}$, so $\Pi(\CL_{m,i},\pi)\not\subseteq H$. Therefore
$$
L_H(\CL_{m,i})\geq
\frac{mk\cdot 2^{k-2} - (k+2)2^{k-3}}{2m-2}
=k\cdot2^{k-3}+\frac{(k-2)2^{k-3}}{2m-2}.
$$
This implies that the exponent $B_k=k\cdot 2^{k-3}$ in Theorem \ref{thm:narrow-kap} is the natural limit of our linear forms conditions.
\end{remark}
 
\section{Proof of Proposition \ref{prop:new}}\label{sec:new}

\subsection{A model case}

First we illustrate the ideas by sketching the proof of (2) in the case $k=3$ and $i=1$. After introducing some short shifts, we can arrange the sum as a bilinear average. As mentioned earlier, we bound this bilinear average by an operator norm and use the inequality $\|T\|_{\operatorname{op}}^{2m} \leq \operatorname{tr}((TT^*)^m)$, instead of the usual procedure of applying Cauchy--Schwarz, which corresponds to using this inequality with $m=2$.

Let $f_2,f_3,g_2,g_3:G \rightarrow\R$ be functions bounded by $\nu_+:=\nu+1$. We need to bound
$$
\E_{n \in G}\E_{d_0,d_1\in [D]} (\nu(n)-1) f_2(n+d_0)f_3(n+2d_0) g_2(n+d_1)g_3(n+2d_1).
$$
Introduce two short shifts $u,v \in [S]$. Shifting $d_0,d_1$ by $u+v$  (which causes an error $O(\eta)$ by the  linear forms conditions) and shifting $n$ by $-u-2v$, we reduce the task to bounding
$$
J =\E_{n,d_0,d_1,u,v}(\nu(n-u-2v)-1) f_2(n+d_0-v)f_3(n+2d_0+u) g_2(n+d_1-v)g_3(n+2d_1+u).
$$
Note that these are the five linear forms \eqref{eq:psik3} in $\Psi_1$:
$$
-u-2v, \ \ d_0-v, \ \ d_1-v, \ \ 2d_0+u, \ \ 2d_1+u.
$$
Since $\Theta=\Theta_1$ is defined to be the set of parameters $\bm{\theta}=(d_0,d_1)$ such that these five linear forms are distinct linear forms in $u,v$, in this model case we just have $\Theta = \{\bm{\theta}=(d_0,d_1): d_0\neq d_1\}$. Before proceeding, let us treat the contribution to $J$ from those $\bm{\theta}\notin\Theta$, which is bounded using the triangle inequality by
\begin{equation}\label{eq:d0=d1}
\E_{n,d_0,d_1,u,v}1_{d_0=d_1} \prod_{\psi\in\Psi_1}\nu_+(n+\psi(d_0,d_1,u,v))
\end{equation}
Expanding each factor $\nu_+=\nu+1$ and using the linear forms conditions for $\Psi_1$, we have
$$
\E_{n,d_0,d_1,u,v}\prod_{\psi\in\Psi_1}\nu_+(n+\psi(d_0,d_1,u,v)) = 2^5+O(\eta).
$$
On the other hand, the five forms in $\Psi_1$ occur as a subcollection of $\CL_{m,1}$ (by taking the two forms at each of the vertices $v_1,u_1$ and the form on the edge joining them; see the figure in Section \ref{sec:linear-forms}). Thus the restricted linear forms conditions for $\CL_{m,1}$ imply that 
$$
\E_{n,d_0,d_1,u,v}1_{d_0\neq d_1}\prod_{\psi\in\Psi_1}\nu_+(n+\psi(d_0,d_1,u,v))  =2^5+O(\eta).
$$
Subtracting the two expressions above, we obtain a bound $O(\eta)$ for the quantity \eqref{eq:d0=d1}. This treats the diagonal contribution from $d_0=d_1$.

We now restrict to $d_0 \neq d_1$.
For $n \in G$ and $\bm{\theta} = (d_0,d_1)\in [D]\times [D]$, define
$$
\begin{aligned}
A_{n,\bm{\theta}}(v) & =f_2(n+d_0-v)g_2(n+d_1-v),\ \ \alpha_{n,\bm{\theta}}(v)=\nu_+(n+d_0-v)\nu_+(n+d_1-v),\\
B_{n,\bm{\theta}}(u)&=f_3(n+2d_0+u)g_3(n+2d_1+u), \ \  \beta_{n,\bm{\theta}}(u)=\nu_+(n+2d_0+u)\nu_+(n+2d_1+u),
\end{aligned}
$$
so that $|A_{n,\bm{\theta}}(v)|\leq\alpha_{n,\bm{\theta}}(v)$ and $|B_{n,\bm{\theta}}(u)|\leq\beta_{n,\bm{\theta}}(u)$.
Define $K_{n,\bm{\theta}}(v,u)=\nu(n-u-2v)-1$. Then
$$
J = \E_{n,\bm{\theta},u,v} K_{n,\bm{\theta}}(v,u) A_{n,\bm{\theta}}(v)B_{n,\bm{\theta}}(u) = \E_{n,\bm{\theta}}1_{\bm{\theta}\in\Theta}J_{n,\bm{\theta}} + O(\eta),
$$
where $J_{n,\bm{\theta}}$ is given by the bilinear average
$$
J_{n,\bm{\theta}} = \E_{u,v}K_{n,\bm{\theta}}(v,u) A_{n,\bm{\theta}}(v)B_{n,\bm{\theta}}(u).
$$
Since $u,v\in [S]$, we form the $S\times S$ matrix $T_{n,\bm{\theta}}$ whose $(v,u)$-entry is
$$
T_{n,\bm{\theta}}(v,u) = S^{-1}\sqrt{\alpha_{n,\bm{\theta}}(v)\beta_{n,\bm{\theta}}(u)}K_{n,\bm{\theta}}(v,u).
$$
Viewing $A_{n,\bm{\theta}},B_{n,\bm{\theta}}$ as column vectors in $\R^S$, we can write
$$
J_{n,\bm{\theta}} = \Big\langle \frac{A_{n,\bm{\theta}}}{\sqrt{\alpha_{n,\bm{\theta}}}}, T_{n,\bm{\theta}} \frac{B_{n,\bm{\theta}}}{\sqrt{\beta_{n,\bm{\theta}}}}\Big\rangle \leq \|T_{n,\bm{\theta}}\|_{\op} (\E_v\alpha_{n,\bm{\theta}}(v))^{1/2} (\E_u\beta_{n,\bm{\theta}}(u))^{1/2},
$$
where $\|T_{n,\bm{\theta}}\|_{\op}$ denotes the operator norm of $T_{n,\bm{\theta}}$ and the inequality follows from Cauchy--Schwarz. Here the inner product is normalized as $\langle\ve{x},\ve{y}\rangle = S^{-1}\sum_v x_vy_v$.
By Cauchy--Schwarz, we have
\begin{equation}\label{eq:J2-bound}
|J|^2 \ll \Big(\E_{n,\bm{\theta}}1_{\bm{\theta}\in\Theta}\|T_{n,\bm{\theta}}\|_{\op}^2\Big) \Big(\E_{n,\bm{\theta}}(\E_v\alpha_{n,\bm{\theta}}(v))(\E_u\beta_{n,\bm{\theta}}(u))\Big) + O(\eta^2).
\end{equation}
The second factor on the right-hand side above is
$$
\E_{n,\bm{\theta},u,v} \nu_+(n+d_0-v)\nu_+(n+d_1-v)\nu_+(n+2d_0+u)\nu_+(n+2d_1+u),
$$
which is $O(1)$ by the linear forms conditions for $\Psi_1$.
To bound the operator norm, we use
$$
\|T_{n,\bm{\theta}}\|_{\mathrm{op}}^{2m}\leq\operatorname{tr}((T_{n,\bm{\theta}} T_{n,\bm{\theta}}^*)^m).
$$
This inequality follows by noting that
$$
\|T_{n,\bm{\theta}}\|_{\op} = s_1, \ \ \operatorname{tr}((T_{n,\bm{\theta}} T_{n,\bm{\theta}}^*)^m)=\sum_js_j^{2m},
$$
where $s_1 \geq s_2 \geq \cdots \geq 0$ are the singular values of $T_{n,\bm{\theta}}$. By H\"{o}lder's inequality, we have
\begin{equation}\label{eq:J-bound}
|J|^{2m} \ll\E_{n,\bm{\theta}} 1_{\bm{\theta}\in\Theta} \|T_{n,\bm{\theta}}\|_{\op}^{2m} + \eta^{2m} \ll 
\E_{n,\bm{\theta}} 1_{\bm{\theta}\in\Theta} \operatorname{tr}((T_{n,\bm{\theta}} T_{n,\bm{\theta}}^*)^m) + \eta^{2m}.
\end{equation}
Expanding the trace, we have
$$
\operatorname{tr}((T_{n,\bm{\theta}} T_{n,\bm{\theta}}^*)^m) =
\E_{\ve u,\ve v\in[S]^m}\prod_{h=1}^m
\alpha_{n,\bm{\theta}}(v_h)\beta_{n,\bm{\theta}}(u_h)
K_{n,\bm{\theta}}(v_h,u_h)K_{n,\bm{\theta}}(v_{h+1},u_h),
$$
with $v_{m+1}=v_1$.
Recalling the definitions of $\alpha_{n,\bm{\theta}},\beta_{n,\bm{\theta}},K_{n,\bm{\theta}}$, for each $h$ we have a product of six factors:
$$
\begin{aligned}
&\nu_+(n+d_0-v_h), \ \ \nu_+(n+d_1-v_h), \ \ \nu_+(n+2d_0+u_h), \ \ \nu_+(n+2d_1+u_h), \\
&(\nu-1)(n-u_h-2v_h), \ \ (\nu-1)(n-u_h-2v_{h+1}).
\end{aligned}
$$
In total there are $6m$ linear forms, which are precisely those in $\CL_{m,1}$ as depicted in the figure in Section \ref{sec:linear-forms}. Thus
$$
\E_{n,\bm{\theta}} 1_{\bm{\theta}\in\Theta} \operatorname{tr}((T_{n,\bm{\theta}} T_{n,\bm{\theta}}^*)^m) = 
\E_{n,\bm{\theta}} \E_{\ve{u},\ve{v}}1_{\bm{\theta}\in\Theta} \prod_{\psi\in\CL_{m,1}}\nu_\psi(n+\psi(\bm{\theta},\ve{u},\ve{v})),
$$
where each $\nu_\psi$ is either $\nu-1$ or $\nu_+$, and at least one of them is $\nu-1$. After expanding and using the restricted linear forms conditions for $\CL_{m,1}$, we get a bound $O(\eta)$ for the expression above. Combining this with \eqref{eq:J-bound}, we complete the proof.

\subsection{The general case: initial applications of Cauchy--Schwarz}\label{sec:new-CS}

We now turn to the general case, which reduces to the same bilinear argument in the model case after $k-3$ applications of Cauchy--Schwarz.
Fix $1 \leq i \leq k$, and let $f_j^{(\tau)}$ be functions on $G$ with $|f_j^{(\tau)}| \leq \nu_+:= \nu+1$ for $1 \leq j \leq k$ and $\tau \in \{0,1\}$. Let 
$$
I = \E_{n \in G}(\nu(n)-1) \CD_i((f_j^{(0)})_{j\neq i})(n) \CD_i((f_j^{(1)})_{j\neq i})(n).
$$
Expanding the definitions of the two dual functions, we have
$$
I = \E_{n \in G}\E_{d_0,d_1\in [D]}(\nu(n)-1) \prod_{j\neq i}\prod_{\tau\in\{0,1\}} f_j^{(\tau)}(n+(j-i)d_\tau).
$$
Let $S = \eta D$ so that $\nu$ satisfies the $(k,m)$-linear forms conditions of level $\eta$ at scale $S$.
Introduce short shifts $\ve{s}=(s_j)_{j\neq i}$ where each $s_j \in[S]$. Shifting both $d_0,d_1$ by $\sum_{j\neq i}s_j$ and shifting $n$ by $\psi_i(\ve{s})=\sum_{\ell\neq i}(i-\ell)s_\ell$, we get
$$
I=\E_{n\in G}\E_{(d_0,d_1,\ve{s})\in\Omega}(\nu(n+\psi_i(\ve s))-1)
\prod_{j\neq i}\prod_{\tau\in\{0,1\}}f_j^{(\tau)}(n+\psi_j^{(\tau)}(d_0,d_1,\ve{s})),
$$
where
$$
\psi_j^{(\tau)}(d_0,d_1,\ve{s})=(j-i)d_\tau+\sum_{l\neq i}(j-l)s_l
$$
and
$$
\Omega = \Big\{(d_0,d_1,(s_j)_{j\neq i}): s_j \in [S], d_0+\sum_{j\neq i}s_j \in [D], d_1+\sum_{j\neq i}s_j \in [D]\Big\}.
$$
This shift ensures that $\psi_j^{(\tau)}$ does not depend on $s_j$ for each $j$. Note that the formula above for $\psi_j^{(\tau)}$ can also be applied to $j=i$ and is consistent with the formula for $\psi_i$, which does not depend on $d_0,d_1$ and so the superscript in $\tau$ is dropped.

We would like to replace the range $\Omega$ for $(d_0,d_1,\ve{s})$ by the rectangular box $[D]^2\times [S]^{k-1}$. Note that if $(d_0,d_1,\ve{s})$ lies in the symmetric difference $\Delta$ between $\Omega$ and the box, then $d_{\tau} \in [-kS,0] \cup [D-kS, D]$ for some $\tau\in\{0,1\}$. Hence $\Delta$ is contained in the union of four boxes, each of which has all side lengths at least $S$ and volume $O(S^kD)$. Since the linear forms $\psi_i$ and $\psi_j^{(\tau)}$ form a subcollection of $\Psi_i$, the linear forms conditions for $\Psi_i$ imply that
$$
\frac{1}{|\Omega|}\E_{n \in G} \sum_{(d_0,d_1,\ve{s})\in\Delta} \nu_+(n+\psi_i(\ve{s}))\prod_{j\neq i}\prod_{\tau\in\{0,1\}}\nu_+(n+\psi_j^{(\tau)}(d_0,d_1,\ve{s})) \ll \frac{|\Delta|}{|\Omega|} \ll \frac{S}{D} \ll \eta.
$$
It follows that we can replace the range $\Omega$ in $I$ by the box at an error of $O(\eta)$:
\begin{equation}\label{eq:I-induction-base}
I=\E_{n,d_0,d_1,\ve s}(\nu(n+\psi_i(\ve s))-1)
\prod_{j\neq i}\prod_{\tau\in\{0,1\}}f_j^{(\tau)}(n+\psi_j^{(\tau)}(d_0,d_1,\ve{s})) + O(\eta),
\end{equation}
where we suppress the ranges of summation for convenience:
$$
n \in G, \ \ d_0,d_1 \in [D], \ \ s_j \in [S]\text{ for }j\neq i.
$$

Recall the choice of $a,b \in [k]\setminus\{i\}$ at the beginning of Section \ref{sec:linear-forms}.
For each index $j\in [k]\setminus\{i,a,b\}$, we will perform an application of Cauchy--Schwarz to eliminate the functions $f_j^{(\tau)}$ (by replacing them with their upper bound $\nu_+$). At the end of this process, only $f_a^{(\tau)}$ and $f_b^{(\tau)}$ remain, as in the model case above.
This procedure is summarized in the following lemma, which is essentially \cite[Lemma 9.1]{Shao2017}. We include a proof to make it clear what linear forms conditions are used.

\begin{lemma}[Gowers-Cauchy--Schwarz]\label{lem:gcs}
Let $E \subset [k]\setminus\{i,a,b\}$ be a subset. We have
$$
|I|^{2^{|E|}} \ll |J_E| + O(\eta^{2^{|E|}}),
$$
where
$$
J_E = \E
\prod_{\bm{\omega}\in\{0,1\}^E}(\nu(n+\psi_i^{(\bm{\omega})})-1) 
 \prod_{\substack{j\in E\\ \tau\in\{0,1\} \\ \bm{\omega}\in\{0,1\}^{E\setminus\{j\}}}}\nu_+(n+\psi_j^{(\tau,\bm{\omega})})
 \prod_{\substack{j\notin E\cup\{i\}\\ \tau\in\{0,1\}\\ \bm{\omega}\in\{0,1\}^E}}f_j^{(\tau)}(n+\psi_j^{(\tau,\bm{\omega})}).
$$
Here the average is over
$$
n \in G, \ \ d_0,d_1 \in [D], \ \ s_j^{(0)},s_j^{(1)}\in [S] \text{ for }j \in E,
\ \ s_j \in [S] \text{ for }j\notin E\cup\{i\}
$$
and the linear forms are given by
$$
\psi_j^{(\tau,\bm{\omega})}  = \psi_j^{(\tau)}(d_0,d_1,\ve{s}^{(\bm{\omega})}),
$$
where $\ve{s}^{(\bm{\omega})}$ for $\bm{\omega}=(\omega_j)_{j\in E}$ denotes the vector whose $j$-th entry is given by $s_j^{(\omega_j)}$ for $j \in E$ and $s_j$ for $j\notin E\cup\{i\}$.
\end{lemma}

Here we are intentionally using the same notation $\psi_j^{(\tau,\bm{\omega})}$ for the linear forms as in \eqref{eq:psi_j}, since they do exactly match up when $E = [k]\setminus\{i,a,b\}$, which is what we will apply the lemma with.

\begin{proof}
We induct on the size $|E|$. When $E$ is empty, the conclusion follows from \eqref{eq:I-induction-base}. Now assume that the conclusion has already been proved for $E$, and we deduce the conclusion for $E \cup \{\ell\}$ for some $\ell\notin E\cup\{i,a,b\}$.
Since the two terms involving $f_\ell^{(0)}$ and $f_\ell^{(1)}$ do not depend on $s_\ell$, we can write $J_E$ as
$$
\E \prod_{\tau,\bm{\omega}}f_\ell^{(\tau)}(n+\psi_\ell^{(\tau,\bm{\omega})}) \cdot\E_{s_\ell} \prod_{\bm{\omega}}(\nu(n+\psi_i^{(\bm{\omega})})-1) \prod_{\substack{j\in E\\ \tau,\bm{\omega}}}\nu_+(n+\psi_j^{(\tau,\bm{\omega})})
\prod_{\substack{j\notin E\cup\{i,\ell\}\\ \tau,\bm{\omega}}} f_j^{(\tau)}(n+\psi_j^{(\tau,\bm{\omega})}),
$$
where the outer average is over $n,d_0,d_1,s_j^{(0)},s_j^{(1)}$ for $j\in E$, and $s_j$ for $j\notin E\cup\{i,\ell\}$. We adopt this convention for all unspecified averages appearing in the remainder of this proof.
By the triangle inequality, we can bound $J_E$ by
$$
\E \prod_{\tau,\bm{\omega}}\nu_+(n+\psi_\ell^{(\tau,\bm{\omega})}) 
\Big|\E_{s_\ell} \prod_{\bm{\omega}}(\nu(n+\psi_i^{(\bm{\omega})})-1) \prod_{\substack{j\in E\\ \tau,\bm{\omega}}}\nu_+(n+\psi_j^{(\tau,\bm{\omega})})
\prod_{\substack{j\notin E\cup\{i,\ell\}\\ \tau,\bm{\omega}}} f_j^{(\tau)}(n+\psi_j^{(\tau,\bm{\omega})})\Big|.
$$
By Cauchy--Schwarz, $|J_E|^2$ is bounded by the product of the two terms
$$
\E \prod_{\tau,\bm{\omega}}\nu_+(n+\psi_\ell^{(\tau,\bm{\omega})})
$$
and
$$
\E \prod_{\tau,\bm{\omega}}\nu_+(n+\psi_\ell^{(\tau,\bm{\omega})}) \Big|\E_{s_\ell} \prod_{\bm{\omega}}(\nu(n+\psi_i^{(\bm{\omega})})-1) \prod_{\substack{j\in E\\ \tau,\bm{\omega}}}\nu_+(n+\psi_j^{(\tau,\bm{\omega})})
\prod_{\substack{j\notin E\cup\{i,\ell\}\\ \tau,\bm{\omega}}} f_j^{(\tau)}(n+\psi_j^{(\tau,\bm{\omega})})\Big|^2.
$$
The first term above is $O(1)$ after expanding and using the linear forms conditions (for an appropriate subcollection of $\Psi_i$). For the second term, we note that $\bm{\omega}$ runs over $\{0,1\}^{E\setminus\{j\}}$ for $j \in E$ and $\bm{\omega}\in\{0,1\}^E$ for other $j$. After expanding the square and naming the two copies of $s_\ell$ as $s_{\ell}^{(0)},s_\ell^{(1)}$, the second term is precisely $J_{E\cup\{\ell\}}$. Thus
$$
|J_E|^2 \ll J_{E\cup\{\ell\}},
$$
and the desired conclusion for $E\cup\{\ell\}$ follows from the induction hypothesis.
\end{proof}

\subsection{The general case: treating the resulting bilinear average}

Let $E = [k]\setminus\{i,a,b\}$ from now on. By Lemma \ref{lem:gcs} applied to $E$, we have
$$
|I|^{2^{k-3}} \ll |J| + O(\eta^{2^{k-3}}),
$$
where
\begin{equation}\label{eq:J}
J=\E
\prod_{\bm{\omega}\in\{0,1\}^E}(\nu(n+\psi_i^{(\bm{\omega})})-1)
\prod_{\substack{j\in E\\ \tau\in\{0,1\} \\ \bm{\omega}\in\{0,1\}^{E\setminus\{j\}}}}\nu_+(n+\psi_j^{(\tau,\bm{\omega})}) \prod_{\substack{j\in\{a,b\}\\ \tau\in\{0,1\} \\ \bm{\omega}\in\{0,1\}^E}}f_j^{(\tau)}(n+\psi_j^{(\tau,\bm{\omega})}).
\end{equation}
Here the average is over 
$$
n \in G, \ \ d_0,d_1 \in [D], \ \ s_a,s_b \in [S], \ \  s_j^{(\tau)} \in [S]\text{ for }j \in E.
$$
As noted earlier, the linear forms appearing in \eqref{eq:J} are precisely those in $\Psi_i$ defined in \eqref{eq:psi_j}, and we continue to adopt the conventions mentioned there. We will treat $J$ as a bilinear average as in the model case, but with a much more complicated operator (matrix) $K_{n,\bm{\theta}}(v,u)$.

Let $\bm{\theta} = (d_0,d_1,(s_j^{(0)},s_j^{(1)})_{j\in E})$, $u =s_a$, and $v=s_b$. Recall that $\Theta$ is the set of parameters $\bm{\theta}$ such that the linear forms in $\Psi_i$ are distinct linear forms in $u,v$.
Before proceeding, we first treat those $\bm{\theta}\notin\Theta$.

\begin{lemma}\label{lem:bad-theta}
Following the notations and assumptions as above, we have
$$
\Big|\E_{n,\bm{\theta},u,v} 1_{\bm{\theta}\notin\Theta_i} \prod_{\psi\in\Psi_i}\nu_+(n+\psi(\bm{\theta},u,v))\Big| = O(\eta).
$$
\end{lemma}

\begin{proof}
By the linear forms conditions for $\Psi_i$, we have
$$
\E_{n,\bm{\theta},u,v}\prod_{\psi\in\Psi_i}\nu_+(n+\psi(\bm{\theta},u,v)) = 2^{|\Psi_i|} + O(\eta)
$$
after expanding $\nu_+=\nu+1$. On the other hand, since $\Psi_i$ is a subcollection of $\CL_{m,i}$, the restricted linear forms conditions for $\CL_{m,i}$ imply that
$$
\E_{n,\bm{\theta},u,v}1_{\bm{\theta}\in\Theta}\prod_{\psi\in\Psi_i}\nu_+(n+\psi(\bm{\theta},u,v)) = 2^{|\Psi_i|} + O(\eta).
$$
The desired conclusion follows by subtracting the two estimates above.
\end{proof}

By Lemma \ref{lem:bad-theta}, the contribution to $J$ in \eqref{eq:J} from $\bm{\theta}\notin\Theta_i$ is $O(\eta)$, and thus we can restrict to $\bm{\theta}\in\Theta_i$:
$$
J = \E_{n,\bm{\theta}}1_{\bm{\theta}\in\Theta_i} J_{n,\bm{\theta}} + O(\eta),
$$
where $J_{n,\bm{\theta}}$ is the bilinear average
$$
J_{n,\bm{\theta}} = \E_{u,v} K_{n,\bm{\theta}}(v,u)A_{n,\bm{\theta}}(v)B_{n,\bm{\theta}}(u),
$$
where
$$
A_{n,\bm{\theta}}(v)=\prod_{\substack{\tau\in\{0,1\} \\ \bm{\omega}\in\{0,1\}^E}}
f_a^{(\tau)}(n+\psi_a^{(\tau,\bm{\omega})}), \ \ 
B_{n,\bm{\theta}}(u)=\prod_{\substack{\tau\in\{0,1\} \\ \bm{\omega}\in\{0,1\}^E}}
f_b^{(\tau)}(n+\psi_b^{(\tau,\bm{\omega})}),
$$
and
$$
K_{n,\bm{\theta}}(v,u)=
\prod_{\bm{\omega}\in\{0,1\}^E}(\nu(n+\psi_i^{(\bm{\omega})})-1)
\prod_{j\in E}\prod_{\substack{\tau\in\{0,1\}\\
\bm{\omega}\in\{0,1\}^{E\setminus\{j\}}}}\nu_+(n+\psi_j^{(\tau,\bm{\omega})}).
$$
From this stage, we follow the argument in the model case, with this more complicated $K_{n,\bm{\theta}}(v,u)$. We have $|A_{n,\bm{\theta}}(v)| \leq \alpha_{n,\bm{\theta}}(v)$ and $|B_{n,\bm{\theta}}(u)| \leq \beta_{n,\bm{\theta}}(u)$, where
$$
\alpha_{n,\bm{\theta}}(v)=\prod_{\substack{\tau\in\{0,1\} \\ \bm{\omega}\in\{0,1\}^E}}
\nu_+(n+\psi_a^{(\tau,\bm{\omega})}), \ \ 
\beta_{n,\bm{\theta}}(u)=\prod_{\substack{\tau\in\{0,1\} \\ \bm{\omega}\in\{0,1\}^E}}
\nu_+(n+\psi_b^{(\tau,\bm{\omega})}).
$$
Since $u,v\in [S]$, we form the $S \times S$ matrix $T_{n,\bm{\theta}}$ whose $(v,u)$-entry is
$$
T_{n,\bm{\theta}}(v,u) = S^{-1}\sqrt{\alpha_{n,\bm{\theta}}(v)\beta_{n,\bm{\theta}}(u)}K_{n,\bm{\theta}}(v,u).
$$
By Cauchy--Schwarz we arrive at \eqref{eq:J2-bound}:
$$
|J|^2 \ll \Big(\E_{n,\bm{\theta}}1_{\bm{\theta}\in\Theta_i}\|T_{n,\bm{\theta}}\|_{\op}^{2}\Big) \Big(\E_{n,\bm{\theta}}(\E_v\alpha_{n,\bm{\theta}}(v))(\E_u\beta_{n,\bm{\theta}}(u))\Big) + O(\eta^2)
$$
The second factor on the right-hand side above is
$$
\E_{n,\bm{\theta},u,v} \prod_{\substack{\tau\in\{0,1\}\\ \bm{\omega}\in\{0,1\}^E}}
\nu_+(n+\psi_a^{(\tau,\bm{\omega})})
\nu_+(n+\psi_b^{(\tau,\bm{\omega})}) = O(1)
$$
by the linear forms conditions for $\Psi_i$.
To bound the first factor involving the operator norm, we use
$$
\|T_{n,\bm{\theta}}\|_{\op}^{2m} \leq \operatorname{tr}((T_{n,\bm{\theta}} T_{n,\bm{\theta}}^*)^m).
$$
By H\"{o}lder's inequality, we arrive at \eqref{eq:J-bound}:
$$
|J|^{2m} \ll \E_{n,\bm{\theta}}1_{\bm{\theta}\in\Theta_i}\|T_{n,\bm{\theta}}\|_{\op}^{2m} + \eta^{2m} \ll \E_{n,\bm{\theta}}1_{\bm{\theta}\in\Theta_i}\operatorname{tr}((T_{n,\bm{\theta}} T_{n,\bm{\theta}}^*)^m) + \eta^{2m}.
$$
The trace can be expanded as
$$
\operatorname{tr}((T_{n,\bm{\theta}} T_{n,\bm{\theta}}^*)^m) =
\E_{\ve u,\ve v\in[S]^m}\prod_{h=1}^m
\alpha_{n,\bm{\theta}}(v_h)\beta_{n,\bm{\theta}}(u_h)
K_{n,\bm{\theta}}(v_h,u_h)K_{n,\bm{\theta}}(v_{h+1},u_h),
$$
with $v_{m+1}=v_1$. Recalling the definitions of $\alpha_{n,\bm{\theta}},\beta_{n,\bm{\theta}},K_{n,\bm{\theta}}$, the linear forms appearing in the above product are precisely those in $\CL_{m,i}$. Thus
$$
\E_{n,\bm{\theta}}1_{\bm{\theta}\in\Theta_i}\operatorname{tr}((T_{n,\bm{\theta}} T_{n,\bm{\theta}}^*)^m) = \E_{n,\bm{\theta},\ve{u},\ve{v}}1_{\bm{\theta}\in\Theta_i}
 \prod_{\psi\in \CL_{m,i}} \nu_{\psi}(n+\psi(\bm{\theta},\ve{u},\ve{v})),
$$
where each $\nu_{\psi}$ is either $\nu-1$ or $\nu_+$, and at least one of them is $\nu-1$. By the restricted linear forms conditions for $\CL_{m,i}$, the quantity above is $O(\eta)$ after expanding each $\nu_{\psi}$. It follows that $|J|^{2m}\ll\eta$ and hence $|I| \ll \eta^{1/(2m\cdot 2^{k-3})}$ as desired.

\subsection{The $L^2$ estimate for dual functions}

We conclude by giving the short proof of Proposition \ref{prop:new}(1).  Fix $1 \leq i \leq k$, and let $\nu_1,\cdots,\nu_k \in \{1,\nu\}$.
Write $F=\CD_i((\nu_j)_{j\neq i})$. It suffices to show that
$$
\E_{n\in G} F(n)=1+O(\eta),\ \  \E_{n\in G} F(n)^2=1+O(\eta).
$$
For the second moment, we have
$$
\E_{n \in G} F(n)^2 = \E_{n,d_0,d_1} \prod_{j\neq i}\prod_{\tau\in\{0,1\}}\nu_j(n+(j-i)d_\tau).
$$
Let $S = \eta D$ so that $\nu$ satisfies the $(k,m)$-linear forms conditions of level $\eta$ at scale $S$. As in the beginning of Section \ref{sec:new-CS}, we introduce short shifts $\ve{s}=(s_j)_{j\neq i}$ where each $s_j \in [S]$, shift both $d_0,d_1$ by $\sum_{j\neq i}s_j$, and shift $n$ by $\sum_{\ell\neq i}(i-\ell)s_{\ell}$. This leads to
$$
\E_{n \in G} F(n)^2 = \E_{n,d_0,d_1,\ve{s}} \prod_{j\neq i}\prod_{\tau\in\{0,1\}}\nu_j(n+\psi_j^{(\tau)}(d_0,d_1,\ve{s})) + O(\eta),
$$
where $O(\eta)$ comes from applying the linear forms conditions for $\Psi_i$ to a symmetric difference as before, and the variables are averaged over the range
$$
n\in G, \ \ d_0,d_1 \in [D], \ \ s_j\in [S]\text{ for }j \neq i.
$$
Since the linear forms $\psi_j^{(\tau)}$ form a subcollection of $\Psi_i$, the linear forms conditions for $\Psi_i$ imply that the average on the right-hand side above is $1+O(\eta)$. It follows that
$$
\E_{n\in G} F(n)^2 = 1 + O(\eta),
$$
as desired. The estimate for the first moment of $F(n)$ is similar (and easier).

\section{The transference argument}\label{sec:transference}

In this section, we give the transference step summarized in Lemmas \ref{lem:dense-model} and \ref{lem:counting}. We follow the arguments in \cite{Zhao2014,Shao2017,ConlonFoxZhao2015}.

\subsection{The dense model argument: Proof of Lemma \ref{lem:dense-model}}

For a positive integer $T$, let $\CF_T$ be the convex hull of the functions $u: G\rightarrow \R$ of the form
\begin{equation}\label{eq:dense-model-u}
u(n) = \E_{\ve{s}=(s_1,\cdots,s_{k-1})\in [T]^{k-1}} \prod_{j=1}^{k-1}u_j\Big(n - k!(s_1+\cdots+s_{k-1}), \ve{s}\Big),
\end{equation}
where each $u_j$ is a $1$-bounded function such that $u_j(n,\ve{s})$ does not depend on $s_j$. For a function $u\in\CF_S$ of the form \eqref{eq:dense-model-u}, we have
$$
|\langle \nu-1,u\rangle|^{2^{k-1}} \leq \E_{n \in G}\E_{\ve{s}^{(0)},\ve{s}^{(1)}\in [S]^{k-1}} \prod_{\bm{\omega}\in\{0,1\}^{k-1}}\Big(\nu\Big(n+k!\sum_{j=1}^{k-1}s_j^{(\omega_j)}\Big)-1\Big).
$$
This follows from $k-1$ repeated applications of Cauchy--Schwarz, with respect to the variable $s_j$ in the $j$th step. Since $u_j$ does not depend on $s_j$, the Cauchy--Schwarz step with respect to $s_j$ eliminates the function $u_j$. Expanding the right-hand side product and using the hypothesis on $\nu$, we have
$$
|\langle \nu-1,u\rangle| \ll \eta^{1/2^{k-1}}.
$$
Thus the same bound holds for all $u \in \CF_S$. Now choose $T = \lfloor \eta^{1/2}D\rfloor$, and we will apply the Green-Tao-Ziegler dense model theorem to $\CF_T$ (see \cite[Lemma 8.5]{Shao2017} or \cite{RTTV2008}). Let $K=K(\eps)$  and $\eps'=\eps'(\eps)$ be the constants from \cite[Lemma 8.5]{Shao2017}. By the approximate algebra property that any product of $K$ functions in $\CF_T$ can be well approximated by a  function in $\CF_S$ (see \cite[Lemma 8.6]{Shao2017}), for each  $u \in \CF_T^K$ there exists $v \in \CF_S$ satisfying $\|v-u\|_{\infty} = O(K\eta^{1/2})$, and thus
$$
|\langle \nu-1,u\rangle| \ll |\langle\nu-1,v\rangle| + \|v-u\|_{\infty} \ll \eta^{1/2^{k-1}} + K\eta^{1/2}.
$$
The bound above can be made at most $\eps'$ since $\eta$ is small enough in terms of $\eps$. Hence the dense model theorem produces a function $\wt{f}:G\rightarrow [0,1]$ with $\E\wt{f}=\E f$, such that
$$
\sup_{u \in \CF_T}|\langle f-\wt{f},u\rangle| \leq \eps.
$$
Fix $i \in [k]$.
In the notation of \cite[Definition 8.1]{Shao2017}, this implies that $(f,\wt{f})$ is an $\eps$-discrepancy pair with width $T$ with respect to the linear form in $k-1$ variables
$$
\xi((s_j)_{j \neq i}) = k!\sum_{j \neq i}s_j.
$$
By \cite[Lemma 8.3]{Shao2017}, $(f,\wt{f})$ is also an $\eps$-discrepancy pair with width $k!\cdot T$ with respect to 
$$
\psi_i((s_j)_{j \neq i}) = \sum_{j \neq i}(i-j)s_j.
$$
Finally, by \cite[Lemma 8.4]{Shao2017} we have $\|f-\wt{f}\|_{D,i} \ll\eps$ as desired.

\subsection{The counting argument: Proof of Lemma \ref{lem:counting}}

Let $\ell$ be the number of $\nu_1,\cdots,\nu_k$ which are equal to $\nu$. We will prove by induction on $\ell$ that 
$$
\Lambda_D(f_1,\cdots,f_k) =\Lambda_D(\wt{f}_1,\cdots,\wt{f}_k) + O(\eta^{1/2^\ell})
$$
for all $f_i,\wt{f}_i$ as in the hypotheses. If $\ell = 0$, then $\nu_1=\cdots=\nu_k=1$ and all $f_i,\wt{f}_i$ are $1$-bounded. In this case, a standard telescoping argument gives
$$
\Big|\Lambda_D(f_1,\cdots,f_k) -\Lambda_D(\wt{f}_1,\cdots,\wt{f}_k)\Big| \ll \sum_{i=1}^k \|f_i-\wt{f}_i\|_{D,i} \ll \eta,
$$ 
as desired. Now suppose that $\ell \geq 1$, and $\nu_i = \nu$ for some $i \in [k]$. Write
$$
F = \CD_i((f_j)_{j \neq i}), \ \ \wt{F} = \CD_i((\wt{f}_j)_{j\neq i}), \ \ G = \CD_i((\nu_j)_{j\neq i}).
$$
Then $0 \leq F \leq G$ and $0 \leq \wt{F} \leq 1$. We have
$$
\Lambda_D(f_1,\cdots,f_k) - \Lambda_D(\wt{f}_1,\cdots,\wt{f}_k) = \langle f_i, F-\wt{F}\rangle + \langle f_i-\wt{f}_i,\wt{F}\rangle.
$$
The second term on the right-hand side above is bounded by $\|f_i-\wt{f}_i\|_{D,i}  \leq \eta$. Thus it suffices to bound $\langle f_i, F-\wt{F}\rangle$. By Cauchy--Schwarz and using $|f_i| \leq \nu$, we have
$$
|\langle f_i,F-\wt{F}\rangle|^2  \leq \Big(\E_{n \in G}\nu(n)\Big)\Big(\E_{n \in G} \nu(n) (F(n)-\wt{F}(n))^2\Big).
$$
The hypotheses on $\nu$ easily imply that $\E\nu \leq 2$ (say). Since
$$
\Big|\E_{n \in G}(\nu(n)-1) (F(n)-\wt{F}(n))^2\Big| \leq 4\eta
$$
by expanding the square and applying hypothesis (2), it suffices to prove that
$$
\E_{n \in G} (F(n)-\wt{F}(n))^2 \ll \eta^{1/2^{\ell-1}}.
$$
For this purpose, let $F^{\sharp}=\min(F,1)$, and use the decomposition
$$
(F-\wt{F})^2 = (F-\wt{F})(F-F^{\sharp}) + (F-\wt{F})(F^\sharp-\wt{F}).
$$
For the first term in the above decomposition, we have
$$
\E (F-\wt{F})(F-F^{\sharp}) \leq \E (G+1)\max(G-1,0) \leq \E |G-1|^2 + 2\E |G-1| \ll \eta
$$
by hypothesis (1). For the second term in the decomposition, we have
$$
\E (F^\sharp-\wt{F})(F-\wt{F}) = \E (F^{\sharp}-\wt{F})\CD_i((f_j)_{j\neq i}) - \E (F^\sharp-\wt{F}) \CD_i((\wt{f}_j)_{j\neq i}),
$$
which can be written as the difference of two counting operators
$$
\Lambda_D(f_1,\cdots,F^\sharp-\wt{F},\cdots,f_k) - \Lambda_D(\wt{f}_1,\cdots,F^\sharp-\wt{F},\cdots,\wt{f}_k),
$$
where $F^{\sharp}-\wt{F}$ appears in the $f_i$-component. Since $F^\sharp-\wt{F}$ is $1$-bounded, we may apply the induction hypothesis (after splitting $F^\sharp-\wt{F}$ into positive and negative parts) to bound the difference by $O(\eta^{1/2^{\ell-1}})$. This concludes the proof.

\section{Proof of Theorems \ref{thm:narrow-kap} and \ref{thm:relative-narrow-kap}}\label{sec:final-deduction}

In this section, we deduce Theorem \ref{thm:relative-narrow-kap} using Lemmas \ref{lem:dense-model}, \ref{lem:counting}, and Proposition \ref{prop:new}, and then deduce Theorem \ref{thm:narrow-kap} using Theorem \ref{thm:relative-narrow-kap} and Proposition \ref{prop:majorant}.

\begin{proof}[Proof of Theorem \ref{thm:relative-narrow-kap}]
Let $\eps > 0$ be small enough in terms of $k,m,\delta$, and assume that $\eta$ is small enough in terms of $\eps$. Since $\E f \leq \E \nu \leq 2$ (say), we may assume that $\E f \leq 1$ by renormalizing (and replacing $\delta$ by $\delta/2$). Let $S = \eta D$, so that $\nu$ satisfies the $(k,m)$-linear forms conditions of level $\eta$ at scale $S$.
The extra mild conditions for $\CC$ in the $(k,m)$-linear forms conditions imply that
$$
\Big|\E_{n \in G}\E_{\ve{s}^{(0)},\ve{s}^{(1)}\in [S]^{k-1}} \prod_{\bm{\omega}\in\{0,1\}^{k-1}}\Big(\nu\Big(n+k!\sum_{j=1}^{k-1}s_j^{(\omega_j)}\Big)-1\Big)\Big| \ll \eta
$$
after expanding out the product. Since $\eta$ is sufficiently small in terms of $\eps$, Lemma \ref{lem:dense-model} implies that there exists $\wt{f}:G\rightarrow [0,1]$ with $\E f = \E \wt{f} \geq \delta$, such that
$$
\|f - \wt{f}\|_{D,i} \leq \eps.
$$
for each $i \in [k]$. By the linear forms conditions and Proposition \ref{prop:new}, we have the following properties:
\begin{enumerate}
\item For $\nu_1,\cdots,\nu_k \in \{1,\nu\}$ and all $1 \leq i \leq k$, we have
$$
\|\CD_i((\nu_j)_{j\neq i})-1\|_2 \ll \eta^{1/2} < \eps
$$
\item For all functions $f_j^{(\tau)}:G\rightarrow\R$ with $|f_j^{(\tau)}| \leq \nu+1$ and all $1 \leq i \leq k$, we have
$$
\Big|\E_{n \in G}(\nu(n)-1) \CD_i((f_j^{(0)})_{j\neq i})(n) \CD_i((f_j^{(1)})_{j\neq i})(n)\Big| \ll \eta^{1/(m\cdot 2^{k-2})} < \eps.
$$
\end{enumerate}
By Lemma \ref{lem:counting}, we have
$$
\Lambda_D(f,\cdots,f) = \Lambda_D(\wt{f},\cdots,\wt{f}) + O(\eps^{1/(m\cdot 2^{2k-2})}).
$$
Since $\Lambda_D(\wt{f},\cdots,\wt{f}) \gg_{k,\delta}1$ by Szemer\'{e}di's theorem and the usual Varnavides averaging argument \cite{Szemeredi1975,Varnavides1959}, the desired conclusion follows since $\eps$ is small enough in terms of $\delta$.
\end{proof}

\begin{proof}[Proof of Theorem \ref{thm:narrow-kap}]
Choose a positive integer $m$ large enough in terms of $k,\eps$. Choose $\eta>0$  small enough in terms of $k,\delta,\eps$. Let $w = w(X)$ be a slowly growing function of $X$, and let $W = \prod_{p \leq w}p \leq (\log X)^{\eps/2}$. We may discard the primes in $A$ below $X^{1/2}$. Choose a reduced residue class $b\Mod{W}$ such that the set
$$
A_{W,b} = \{n \in \Z: Wn+b \in A, Wn+b > X^{1/2}\}
$$
satisfies
$$
|A_{W,b}| \geq \frac{\delta X}{2\varphi(W)\log X}.
$$
Choose a prime $N \in [4X/W, 8X/W]$, and let $G = \Z/N\Z$. By Proposition \ref{prop:majorant}, there exists a non-negative function $\nu=\nu_{W,b}: G\rightarrow [0,+\infty)$ with the following properties:
\begin{enumerate}
\item (Majorant property) For $n \in [N]$, if $Wn+b$ is prime and $Wn+b>N^{1/2}$ then
$$
\nu(n) \geq c_0 \frac{\varphi(W)}{W}\log N
$$
for some constant $c_0=c_0(k,m) > 0$.
\item (Linear forms conditions) The majorant $\nu$ satisfies the $(k,m)$-linear forms conditions of level $\eta$ at some scale
$$
S \ll (\log N)^{B_k+\eps/2}.
$$
\end{enumerate}

Let $D = S/\eta$ and we may ensure that $D \geq \eta^{-2}$.
For $n \in [N]$, define
$$
f(n) = \begin{cases} c_0 \frac{\varphi(W)}{W}\log N & \text{if }Wn+b\in A\text{ and }Wn+b>X^{1/2},\\ 0 & \text{otherwise.}\end{cases}
$$
Then $0 \leq f \leq \nu$ and $\E f \gg\delta$. By the relative Szemer\'{e}di theorem (Theorem \ref{thm:relative-narrow-kap}), we have
$$
\Lambda_D(f,\cdots,f) \gg 1.
$$
Hence there exists a non-trivial $k$-term arithmetic progression $n, n+d, \cdots,n+(k-1)d$ in $G$ with $d \in [D]$, such that $f(n+jd) > 0$ for each $j$. Since $D = o(N)$ and $f$ is supported on $[N/4]$, it is an arithmetic progression in the integers. This leads to a desired $k$-term progression
$$
Wn+b, Wn+b+Wd, \cdots, Wn+b+(k-1)Wd
$$
in primes, whose common difference is at most $WD \leq WS/\eta \ll (\log N)^{B_k+\eps}$.
\end{proof}

\bibliographystyle{plain}
\bibliography{narrow-kap-references}

\end{document}